\documentclass[12pt, twoside]{article}
\usepackage{times}
\usepackage{graphics}
\usepackage{theorem}
\usepackage{graphicx}
\usepackage{amsmath}
\usepackage{latexsym}
\usepackage{amssymb,mathrsfs}
\usepackage[flushmargin]{footmisc}

\theorembodyfont{\itshape}
\newtheorem{thm}{Theorem}[section]
\newtheorem{lem}{Lemma}[section]
\newtheorem{prop}{Proposition}[section]

\theorembodyfont{\rmfamily}
\newtheorem{defn}{Definition}[section]{\bf}{\rm}
\newtheorem{assumpt}{Assumption}[section]{\bf}{\rm}

\newtheorem{rem}{Remark}[section]{\itshape}{\rmfamily}
\newenvironment{proof}{\noindent{\it Proof.~~}}{\qed}

\makeatletter
\def\eqnarray{\stepcounter{equation}\let\@currentlabel=\theequation
\global\@eqnswtrue
\global\@eqcnt\z@\tabskip\@centering\let\\=\@eqncr
$$\halign to \displaywidth\bgroup\@eqnsel\hskip\@centering
  $\displaystyle\tabskip\z@{##}$&\global\@eqcnt\@ne 
  \hfil$\;{##}\;$\hfil
  &\global\@eqcnt\tw@ $\displaystyle\tabskip\z@{##}$\hfil 
   \tabskip\@centering&\llap{##}\tabskip\z@\cr}
\makeatother
\makeatletter
    \renewcommand{\theequation}{%
    \thesection.\arabic{equation}}
    \@addtoreset{equation}{section}
  \makeatother

\def\vc#1{\mbox{\boldmath $#1$}}
\def\svc#1{\mbox{\boldmath $\scriptstyle #1$}}

\newcommand{\down}[2]{\smash{\lower#2\hbox{#1}}}
\newcommand{\up}[2]{\smash{\lower-#2\hbox{#1}}}
\newcommand{\dm}{\displaystyle}

\newcommand{\qed}{\hspace*{\fill}$\Box$\rule[-10pt]{0pt}{10pt}}

\newcommand{\EE}{\mathsf{E}}
\newcommand{\PP}{\mathsf{P}}

\newcommand{\calH}{\mathcal{H}}
\newcommand{\calL}{\mathcal{L}}

\newcommand{\calR}{\mathcal{R}}
\newcommand{\calS}{\mathcal{S}}

\newcommand{\bbL}{\mathbb{L}}

\newcommand{\bbM}{\mathbb{M}}
\newcommand{\bbN}{\mathbb{N}}

\newcommand{\bbZ}{\mathbb{Z}}

\newcommand{\rmt}{{\rm t}}
\newcommand{\rmd}{{\rm d}}
\newcommand{\rme}{{\rm e}}

\def\simhm#1{\stackrel{#1}{\sim}}

\def\ooverline#1{\overline{\overline{#1}}}

\newcommand{\dd}[1]{\if#11 1\!\!1 
\else {\if#1C I\!\!\!C
\else {\if#1G I\!\!\!G 
\else {\if#1J J\!\!\!J 
\else {\if#1S S\!\!\!S
\else {\if#1Z Z\!\!\!Z
\else {\if#1Q O\!\!\!\!Q
\else I\!\!#1
\fi} 
\fi}
\fi}
\fi} 
\fi} 
\fi} 
\fi} 

\makeatother

\begin{document}\thispagestyle{plain} 

\hfill
{\small Last update date: \today}

{\Large{\bf
\begin{center}
Subexponential tail equivalence of the stationary queue length distributions of 
BMAP/GI/1 queues with and without retrials
\end{center}
}
}

\begin{center}
{
Hiroyuki Masuyama%
\footnote[2]{E-mail: masuyama@sys.i.kyoto-u.ac.jp}
}

\medskip

{\small
Department of Systems
Science, Graduate School of Informatics, Kyoto University\\
Kyoto 606-8501, Japan
}

\bigskip
\medskip

{\small
\textbf{Abstract}

\medskip

\begin{tabular}{p{0.85\textwidth}}
The main contribution of this paper is to prove the
subexponential tail equivalence of the stationary queue length
distributions in the BMAP/GI/1 queues with and without retrials. We first present a stochastic-decomposition-like
result of the stationary queue length in the BMAP/GI/1 retrial queue,
which is an extension of the stochastic decomposition of the
stationary queue length in the M${}^X$/GI/1 retrial queue.  The
stochastic-decomposition-like result shows that the stationary queue
length distribution in the BMAP/GI/1 retrial queue is decomposed into
two parts: the stationary conditional queue length distribution given
that the server is idle; and a certain matrix sequence associated with
the stationary queue length distribution in the corresponding standard
BMAP/GI/1 queue (without retrials). Using the
stochastic-decomposition-like result and matrix analytic methods, we
prove the subexponential tail equivalence of the stationary queue
length distributions in the BMAP/GI/1 queues with and without
retrials. This tail equivalence result does not necessarily require
that the size of an arriving batch is light-tailed, unlike Yamamuro's
result for the M${}^X$/GI/1 retrial queue (Queueing
Syst.\ 70:187--205, 2012). As a by-product, the key lemma to the proof of the main theorem presents a subexponential asymptotic formula for the stationary distribution of a level-dependent M/G/1-type Markov chain, which is the first reported result on the subexponential asymptotics of level-dependent block-structured Markov chains.
\end{tabular}
}
\end{center}

\begin{center}
\begin{tabular}{p{0.90\textwidth}}
{\small

\medskip

}
\end{tabular}

\end{center}

\section{Introduction}\label{introduction}

Retrial queues are queueing models such that a customer finding all
the servers busy on arrival joins the virtual waiting line (called
{\it orbit}) and retries to occupy an idle server after a random time
(this process is repeated until the customer finds an idle server and
occupies it). Many researchers have studied retrial queues for more
than a half of century since the early studies, e.g.,
\cite{Cohe57,Kost47}. However, the analytical results of retrial
queues are less extensive than those of standard (work-conserving and
non-preemptive) queueing models without retrials. In particular, exact
(that is, not approximate) solutions have been derived for a few
simple models such as M/M/$c$ ($c=1,2,3,4$) retrial queues (see
\cite{Fali97,Hans87,Pear89,Tuan09}). For detailed overview, see the
survey papers \cite{Fali90,Yang87} and the books \cite{Arta08,Fali97}.
See also the bibliographies \cite{Arta99,Arta10} and the references
therein.

Recently, the asymptotic analysis has been a hot topic in the study of
retrial queues. Liu and Zhao~\cite{Liu10} and Kim et al.~\cite{Kim12a}
study the light-tailed asymptotics of the stationary queue length
distribution in the M/M/$c$ retrial queue. These results are extended
to an M/M/$c$ retrial queue with non-persistent customers 
\cite{Kim12b,Liu12}. Kim et al.~\cite{Kim07} study the light-tailed
asymptotics of the stationary queue length distribution in an M/GI/1
retrial queue with exponential retrials, which is generalized to the
Markovian arrival case by Kim et al.~\cite{Kim10b}.

As for the subexponential asymptotics, there are a few studies.
Before reviewing them, we introduce the subexponential class of
distributions and related ones.
\begin{defn}\label{defn-subexp}
\hfill
\begin{enumerate}
\item The nonnegative random variable $U$ and its distribution $F_U$
  are said to be heavy-tailed (denoted by $U \in \calH$ and $F_U \in
  \calH$) if $\lim_{x\to\infty}\rme^{\varepsilon x}\PP(U>x) = \infty$
  for all $\varepsilon > 0$.
\item The nonnegative random variable $U$ and its distribution $F_U$
  are said to be long-tailed (denoted by $U \in \calL$ and $F_U \in
  \calL$) if $\PP(U>x) > 0$ for all $x \ge 0$ and
\[
\lim_{x\to\infty}{\PP(U > x+y)  \over \PP(U>x)} = 1
\quad \mbox{for some (thus all) $y > 0$}.
\]
\item The nonnegative random variable $U$ and its distribution $F_U$
  are said to be subexponential (denoted by $U \in \calS$ and $F_U \in
  \calS$) if $\PP(U>x) > 0$ for all $x \ge 0$ and
\[
\lim_{x\to\infty}{\PP(U_1 + U_2 > x)  \over \PP(U>x)} = 2,
\]
where $U_i$'s ($i=1,2,\dots$) are independent copies of $U$.
\item The nonnegative random variable $U$ and its distribution $F_U$
  belong to class $\calR(-\alpha)$ ($\alpha \ge 0$) if $\PP(U>x)$ is
  regularly varying with index $-\alpha$, i.e.,
\[
\lim_{x\to\infty}
{ \PP(U > vx) \over \PP(U>x)} = v^{-\alpha} \quad \mbox{for all $ v > 0$}.
\]
\end{enumerate}
\end{defn}

It is known that $\cup_{\alpha\ge0}\calR(-\alpha) \subset \calS
\subset \calL \subset \calH$. In particular, class $\calS$ is the
largest tractable subclass of heavy-tailed distributions, and it
includes heavy-tailed Weibull, lognormal, Burr, loggamma
distributions, and Pareto distributions, etc. For further details, see
\cite{Foss11,Gold98,Sigm99}.

We now review the literature on the subexponential asymptotics of
retrial queues.  Kim et al.~\cite{Kim10a} consider an M/GI/1 retrial
queue with exponential retrials and the service time distribution in
$\calR(-\beta)$, where $\beta > 1$. For this retrial queue, the
authors show that the waiting time distribution belongs to class
$\calR(-\beta+1)$. Shang et al.~\cite{Shan06} also consider the M/GI/1
retrial queue with exponential retrials, and they prove the
subexponential tail equivalence of the stationary queue length
distributions in the M/GI/1 queues with and without retrials. In order
to specify this tail equivalence result, we denote by $L^{(\mu)}$ the
stationary queue length in the M/GI/1 retrial queue with exponential
retrial rate $\mu$, and denote by $L^{(\infty)}$ the stationary queue
length in the corresponding standard M/GI/1 queue (it is shown that
$L^{(\mu)}$ converges to $L^{(\infty)}$ in distribution as $\mu \to
\infty$; see Theorem~1.8 in \cite{Fali97}). In this setting, Shang et
al.~\cite{Shan06}'s result is stated as follows: If $L^{(\infty)}
\in\calS$, then
\begin{equation}
\PP(L^{(\mu)} > x) \simhm{x} \PP(L^{(\infty)} > x),
\label{tail-equivalence}
\end{equation}
where $f(x) \simhm{x} g(x)$ represents $\lim_{x\to\infty}f(x)/g(x) =
1$.  Note here that $L^{(\infty)} \in \calS$ and
(\ref{tail-equivalence}) imply $L^{(\mu)} \in \calS$ (see, e.g.,
\cite[Proposition 2.8]{Sigm99}).  Yamamuro~\cite{Yama12} extends the
tail equivalence (\ref{tail-equivalence}) to the batch arrival model,
i.e., M${}^X$/GI/1 retrial queue with exponential retrials, though the
batch size distribution is assumed to be light-tailed.

This paper considers a BMAP/GI/1 retrial queue with exponential
retrials, where BMAP represents batch Markovian arrival process
\cite{Luca91}. The main contribution of this paper is to prove the
subexponential tail equivalence of the stationary queue length
distributions in the BMAP/GI/1 queues with and without retrials, which
is an extension of Yamamuro~\cite{Yama12}'s result.

To prove the main result of this paper, we first present a
stochastic-decomposition-like result of the stationary queue length,
which is a generalization of the stochastic decomposition for the
M${}^X$/GI/1 retrial queue with exponential retrials \cite{Yama12}.
The stochastic-decomposition-like result shows that the stationary
queue length distribution in a BMAP/GI/1 retrial queue with
exponential retrials is decomposed into two parts. The first part is
the stationary conditional queue length distribution given that the
server is idle. On the other hand, the second part itself does not
have a probabilistic interpretation. However, pre-multiplying the
second part by a certain probability vector, we have the stationary
queue length distribution in the corresponding standard BMAP/GI/1
queue (without retrials).

Next we prove the main theorem on the subexponential tail equivalence
by combining the stochastic-decomposition-like result with matrix
analytic methods \cite{He14,Lato99,Neut89}. The key to the proof of
the main theorem is to discuss the tail asymptotics of the stationary
conditional queue length distribution given that the server is idle,
which is reduced, by change of measure, to the subexponential
asymptotics of a level-dependent M/G/1-type Markov chain with
asymptotic level-independence.  To the best of our knowledge, there
are no studies on the subexponential asymptotics of level-dependent
block-structured Markov chains. In
addition, the main theorem of this paper is proved without
Yamamuro~\cite{Yama12}'s assumption mentioned above, i.e., the
light-tailedness of the batch size distribution.

The rest of this paper is divided into three sections. Section~\ref
{sec-preliminary} introduces basic definitions, notation and
preliminary results. Section~\ref{sec-main-result} presents the main
theorem. Section~\ref{proof-lem-lim-overline{p}_0(k)} is devoted to
the proof of a lemma, which is key to prove the main theorem.

\section{Preliminary}\label{sec-preliminary}

\subsection{Basic definitions and notation}

Let $\bbZ = \{0,\pm 1, \pm 2,\dots\}$ and $\bbZ_+ = \{0,1,2,\dots\}$,
$\bbN = \{1,2,3,\dots\}$. 

Let $\vc{e}$ and $\vc{I}$ denote the column vector of ones and the
identity matrix, respectively, with appropriate dimensions according
to the context.  The superscript ``$\rmt$" represents the transpose
operator for vectors and matrices.  The notation $[\,\cdot\,]_{i,j}$
(resp.\ $[\,\cdot\,]_i$) represents the $(i,j)$th (resp.\ $i$th)
element of the matrix (resp.\ vector) in the square brackets.

For any matrix $\vc{M}$, let $|\vc{M}|$ denote a matrix obtained by
taking the absolute value of each element of $\vc{M}$, i.e.,
$[|\vc{M}|]_{i,j} = |[\vc{M}]_{i,j}|$. For any matrix sequence
$\{\vc{M}(k);k\in\bbZ\}$, we define $\overline{\vc{M}}(k)$ and
$\ooverline{\vc{M}}(k)$ ($k \in \bbZ$) as
\[
\overline{\vc{M}}(k) = \sum_{l=k+1}^{\infty}\vc{M}(l),
\qquad
\ooverline{\vc{M}}(k) = \sum_{l=k+1}^{\infty}\overline{\vc{M}}(l),
\]
respectively. We then define the convolution of two matrix sequences
$\{\vc{M}(k);k\in\bbZ\}$ and $\{\vc{N}(k);k\in\bbZ\}$ as follows:
\[
\vc{M} \ast \vc{N}(k) =
\sum_{l\in\bbZ} \vc{M}(k-l) \vc{N}(l),\qquad k\in\bbZ,
\]
where the product $\vc{M}(k_1)\vc{N}(k_2)$ is well-defined for all
$k_1, k_2 \in \bbZ$. For any square matrix sequence
$\{\vc{M}(k);k\in\bbZ\}$, we also define $\{\vc{M}^{\ast
  m}(k);k\in\bbZ\}$ ($m\in\bbN$) as the $m$-fold convolution of
$\{\vc{M}(k)\}$ with itself, i.e.,
\[
\vc{M}^{\ast m}(k) = \sum_{l\in\bbZ}
\vc{M}^{\ast (m-1)}(k-l) \vc{M}(l),
\qquad k \in \bbZ,
\]
where $\vc{M}^{\ast0}(0) = \vc{I}$ and $\vc{M}^{\ast0}(k) = \vc{O}$
for $k \in \bbZ\setminus\{0\}$. In addition, for two matrix-valued
functions $\vc{M}_1(\,\cdot\,)$ and $\vc{M}_2(\,\cdot\,)$ with the
same dimension, the notation $\vc{M}_1(x) \simhm{x} \vc{M}_2(x)$
represents $[\vc{M}_1(x)]_{i,j} \simhm{x} [\vc{M}_2(x)]_{i,j}$, i.e.,
\[
\lim_{x\to\infty}
{ [ \vc{M}_1(x) ]_{i,j} \over [\vc{M}_2(x)]_{i,j} }
 = 1\quad \mbox{for all $i$'s and $j$'s}.
\]

The above definitions and notation for matrices are applied to vectors
and scalars in an appropriate manner.

\subsection{Subexponential asymptotics for BMAP/GI/1 queue without retrials}\label{subsec-BMAP/GI/1}

We first introduce the BMAP. Behind the BAMP, there exists a
continuous-time Markov chain with a finite state space
$\bbM:=\{1,2,\dots,M\}$, which is called the background Markov chain
(or the underlying Markov chain).  Let $\{J(t);t \ge 0\}$ denote the
background Markov chain of the BMAP. Let $N(t)$ ($t \ge 0$) denote the
total number of arrivals in time interval $(0,t]$, where $N(0) = 0$ is
  assumed.

For simplicity, we denote by $E_{\downarrow0}$ an appropriate
real-valued function on $[0,\infty)$ such that
  $\lim_{x\downarrow0}E_{\downarrow0}(x)/x=0$.  It then follows by
  definition (see, e.g., \cite{Luca91}) that the stochastic process
  $\{(N(t),J(t));t\ge0\}$ evolves as follows:
\begin{eqnarray}
\lefteqn{
\PP(
N(t+\Delta t) - N(t) = k, J(t) =j \mid J(0) = i
)
}
\qquad &&
\nonumber
\\
&=& 
\left\{
\begin{array}{ll}
1 + [\vc{C}]_{i,i}\Delta t + E_{\downarrow0}(\Delta t), & k = 0,~i=j \in \bbM,
\\
{}
[\vc{C}]_{i,j} \Delta t + E_{\downarrow0}(\Delta t), & k = 0,~i \neq j,~ i,j \in \bbM,
\\
{}
[\vc{D}(k) ]_{i,j}\Delta t + E_{\downarrow0}(\Delta t), & k \in \bbN,~i,j \in \bbM,
\end{array}
\right.
\label{defn-BMAP}
\end{eqnarray}
where $\vc{D}(k)$ ($k\in\bbN$) is an $M \times M$ nonnegative matrix
and $\vc{C}$ is an $M \times M$ matrix such that $[\vc{C}]_{i,i} < 0$
($i \in \bbM$), $[\vc{C}]_{i,j} \ge 0$ ($i \neq j,~ i,j \in \bbM$) and
$(\vc{C} + \sum_{k=1}^{\infty} \vc{D}(k))\vc{e}=\vc{0}$. The BMAP
characterized in (\ref{defn-BMAP}) is denoted by BMAP
$\{\vc{C},\vc{D}(k);k\in\bbN\}$.

Let $\widehat{\vc{D}}(z) = \sum_{k=1}^{\infty}z^k \vc{D}(k)$. From
(\ref{defn-BMAP}), we then have
\begin{equation}
\EE\left[
z^{N(t)} \cdot \dd{1}(J(t) =j) \mid J(0) = i
\right]
=
\left[
\rme^{(\svc{C} + \widehat{\svc{D}}(z)) t}
\right]_{i,j},
\label{eqn-counting-BAMP}
\end{equation}
where $\dd{1}(\,\cdot\,)$ denotes an indicator function that takes
value of one if the statement in the parentheses is true; and takes
value of zero otherwise.  Note here that the infinitesimal generator
of the background Markov chain $\{J(t);t\ge0\}$ is given by $\vc{C} +
\vc{D}$, where $\vc{D} = \sum_{k=1}^{\infty} \vc{D}(k)$.

We assume that $\vc{C} + \vc{D}$ is irreducible, and
then define $\vc{\pi} > \vc{0}$ as the unique stationary probability
vector of $\vc{C} + \vc{D}$.  We also define $\lambda$ as the mean
arrival rate, i.e.,
\begin{equation}
\lambda 
= \vc{\pi} \sum_{k=1}^{\infty} k\vc{D}(k)\vc{e}
= \vc{\pi} \sum_{k=0}^{\infty} \overline{\vc{D}}(k)\vc{e}.
\label{defn-lambda}
\end{equation}
To exclude trivial cases, we assume $\lambda > 0$, which implies that
\begin{equation}
\vc{D}(k_0) \ge \vc{O}, \neq \vc{O} \quad \mbox{for some $k_0 \in \bbN$}.
\label{assumpt-lambda>0}
\end{equation}

Next we describe a standard BMAP/GI/1 queue, i.e, BMAP/GI/1 queue
without retrials. The system has a single server and a buffer of
infinite capacity. Customers arrive at the system according to BMAP
$\{\vc{C},\vc{D}(k);k\in\bbN\}$.  If customers arriving in a batch
find the server idle, then one of them immediately occupies the server
and the others join the waiting line; otherwise all of them join the
waiting line. We assume that the service times of customers are
independent of BMAP $\{\vc{C},\vc{D}(k);k\in\bbN\}$ and independent
and identically distributed (i.i.d.) according to a general
distribution function $H$ on $[0,\infty)$ with mean $h \in
  (0,\infty)$.

We define $\rho$ as the traffic intensity, i.e., 
\begin{equation}
\rho = \lambda h.
\label{defn-rho}
\end{equation}
We also define $\vc{A}(k)$ ($k \in \bbZ_+$) as an $M \times M$ matrix
such that 
\[
[\vc{A}(k)]_{i,j} = \PP(N(T) = k, J(T) = j \mid J(0) = i),
\qquad i,j\in\bbM,
\]
where $T$ denotes a generic random variable for i.i.d.\ service times
with distribution function $H$. It follows from
(\ref{eqn-counting-BAMP}) that
\begin{equation}
\widehat{\vc{A}}(z) :=
\sum_{k=0}^{\infty}z^k \vc{A}(k)=
\int_0^{\infty} 
\rme^{(\svc{C} + \widehat{\svc{D}}(z)) x}\rmd H(x).
\label{defn-A(k)}
\end{equation}
Note here that $\vc{A}:= \widehat{\vc{A}}(1) = \int_0^{\infty}
\rme^{(\svc{C} + \svc{D}) x}\rmd H(x) > \vc{O}$ and
$\vc{A}\vc{e}=\vc{e}$ because $\vc{C} +\vc{D}$ is an irreducible
infinitesimal generator. Thus $\vc{A}$ has the unique stationary
probability vector, which is equal to $\vc{\pi}$. Note also that
\begin{equation}
\rho = \vc{\pi}\sum_{k=1}^{\infty} k \vc{A}(k) \vc{e}.
\label{eqn-rho}
\end{equation}

Throughout this paper, we assume $\rho < 1$, which is the stability
condition for the standard BMAP/GI/1 queue (see \cite{Loyn62}). We
then summarize the results on the subexponential asymptotics of the
stationary queue length distribution in the standard BMAP/GI/1 queue.

Let $\vc{x}(k)$ ($k\in\bbZ_+$) denote a $1 \times M$ vector whose
$i$th element $[\vc{x}(k)]_i$ ($i\in\bbM$) represents the stationary
joint probability that the queue length in the standard BMAP/GI/1
queue is equal to $k$ and the background Markov chain is in state
$i$. According to \cite{Taki00}, $\{\vc{x}(k);k\in\bbZ_+\}$ is
equivalent to the stationary distribution of the following M/G/1-type
Markov chain:
\[
\vc{P}_{\rm M/G/1}
:=
\left(
\begin{array}{ccccc}
\vc{A}(0) & \vc{A}(1) & \vc{A}(2) & \vc{A}(3) & \cdots
\\
\vc{A}(0) & \vc{A}(1) & \vc{A}(2) & \vc{A}(3) & \cdots
\\
\vc{O}    & \vc{A}(0) & \vc{A}(1) & \vc{A}(2) & \cdots
\\
\vc{O}    & \vc{O}    & \vc{A}(0) & \vc{A}(1) & \cdots
\\
\vdots    & \vdots    &  \vdots   & \vdots   & \ddots
\end{array}
\right).
\]

To characterize $\{\vc{x}(k)\}$, we introduce some matrices.  Let
$\vc{G}$ denote the minimal nonnegative solution of
\[
\vc{G} = \sum_{k=0}^{\infty}\vc{A}(k)\vc{G}^k.
\]
Since $\vc{A}$ is irreducible, $\vc{G}$ is stochastic under the
stability condition $\rho < 1$ (see Theorem 2.3.1 in
\cite{Neut89}). We can also show that $\vc{G} > \vc{O}$ (see page 382
of \cite{Masu09}). Thus $\vc{G}$ has the unique and positive
stationary probability vector, denoted by $\vc{g} > \vc{0}$. Further
it is known \cite{Taki00} that
\begin{equation}
\widehat{\vc{x}}(z)
= (1-\rho)\vc{g}(z-1)
\left( z\vc{I} - \widehat{\vc{A}}(z) \right)^{-1}\widehat{\vc{A}}(z),
\quad |z| < 1.
\label{eqn-hat{x}^{(infty)}(z)}
\end{equation}

\begin{rem}\label{rem-G}
Since $\vc{G} > \vc{O}$, $\vc{G}$ is an irreducible and aperiodic
stochastic matrix. Thus $\vc{G}$ has a simple eigenvalue $\gamma_1=1$
(called Perron-Frobenius eigenvalue) and the other eigenvalues
$\gamma_2, \gamma_3,\cdots,\gamma_M$ are located inside the unit
circle in the complex plane, i.e., $1 > |\gamma_2| \ge |\gamma_3| \ge
\cdots \ge |\gamma_M|$ (see, e.g., \cite[Theorem 8.4.4]{Horn90}). In
addition, $\lim_{m\to\infty}\vc{G}^m = \vc{e}\vc{g} > \vc{O}$ (see,
e.g., \cite[Theorem 8.2.8]{Horn90}).
\end{rem}

Let $\vc{R}(0) = \vc{O}$ and $\vc{R}(k)$ $(k\in\bbN)$ denote
\begin{align*}
\vc{R}(k)
&=
\dm\sum_{m=0}^{\infty}
\vc{A}(k+m+1)\vc{G}^m(\vc{I}-\vc{U}(0))^{-1}, & k &\in \bbN,
\end{align*}
where
\[
\vc{U}(0) = \sum_{m=0}^{\infty}
\vc{A}(m+1)\vc{G}^m.
\]
The matrices $\vc{R}(k)$, $\vc{G}$ and $\vc{U}(0)$ satisfy the
following equation (called {\it $RG$-factorization}; see
\cite[Theorem 14]{Zhao03}):
\begin{equation}
z\vc{I} - \widehat{\vc{A}}(z) 
= \left( \vc{I} - \widehat{\vc{R}}(z) \right) 
(\vc{I}-\vc{U}(0)) (z\vc{I} - \vc{G}),
\label{RG-factorization}
\end{equation}
where $\widehat{\vc{R}}(z)=\sum_{k=1}^{\infty}z^k\vc{R}(k)$.
From (\ref{RG-factorization}), we have
\[
\vc{\pi} 
= (1-\rho)\vc{g} (\vc{I} - \vc{U}(0))^{-1} (\vc{I} - \vc{R})^{-1},
\]
where $\vc{R} = \widehat{\vc{R}}(1)$ (see \cite[Lemma
  14]{Taki03}). Further substituting (\ref{RG-factorization}) into
(\ref{eqn-hat{x}^{(infty)}(z)}) yields
\begin{equation}
\widehat{\vc{x}}(z)
= (1-\rho)\vc{g} \left(\vc{I} - \vc{U}(0) \right)^{-1}
\left( \vc{I} - \widehat{\vc{R}}(z) \right)^{-1}
\widehat{\vc{A}}(z).
\label{eqn-hat{x}^{(infty)}(z)-02}
\end{equation}

We now make the following assumption:
\begin{assumpt}\label{assumpt-subexp}
There exists some $\bbZ_+$-valued random variable $Y \in \calS$ such
that for some nonzero vector $\vc{c}^{\rm A}$,
\begin{equation}
\lim_{k\to\infty}{\ooverline{\vc{A}}(k)\vc{e} \over \PP(Y > k)} 
= \vc{c}^{\rm A} \ge \vc{0},\neq\vc{0}.
\label{cond-prop-Masu13-ANNOR}
\end{equation}
\end{assumpt}

Under Assumption~\ref{assumpt-subexp}, we have the following
result.
\begin{prop}[Corollary 4.1 in \cite{Masu13-ANNOR}]\label{prop-Masu13-ANNOR}
If Assumption~\ref{assumpt-subexp} holds, then
\[
\lim_{k\to\infty}{\overline{\vc{x}}(k) \over \PP(Y > k)} 
= {\vc{\pi}\vc{c}^{\rm A} \over 1 - \rho}\cdot \vc{\pi}.
\]
\end{prop}

Finally we present a sufficient condition
Assumption~\ref{assumpt-subexp}.  To this end, we define $T_{\rme}$ as
a random variable that is independent of BMAP
$\{\vc{C},\vc{D}(k);k\in\bbN\}$ and is distributed with
\[
\PP(T_{\rme} \le x)
= { 1 \over \EE[T] }\int_0^x \PP(T > y) \rmd y
= { 1 \over h }\int_0^x \overline{H}(y) \rmd y,
\quad x \ge 0,
\]
which is called {\it the residual service time} or {\it the
  equilibrium random variable of the service time $T$}.
\begin{prop}\label{lem-assumpt-subexp}
Suppose that (i) $\sqrt{T_{\rme}} \in \calL$; and (ii)
$\lim_{k\to\infty}\rme^{\delta \sqrt{k}}\overline{\vc{D}}(k) < \infty$
for some $\delta > 0$. If $T_{\rme} \in \calS$, then
Assumption~\ref{assumpt-subexp} holds for $Y = \lambda T_{\rme}$ and
$\vc{c}^{\rm A} = \rho\vc{e}$.
\end{prop}

\begin{rem}
According to Proposition~\ref{lem-assumpt-subexp},
Assumption~\ref{assumpt-subexp} does not necessarily require that
$\{\vc{D}(k)\}$ is light-tailed.  Some other sufficient conditions for
Assumption~\ref{assumpt-subexp} are presented in Section 4 in
\cite{Masu13-ANNOR}.
\end{rem}

\noindent
{\it Proof of Proposition~\ref{lem-assumpt-subexp}}.~ Let
$\vc{A}_{\rme}(k)$ ($k\in\bbZ_+$) denote an $M \times M$ matrix such
that
\[
[\vc{A}_{\rme}(k)]_{i,j}
= \PP(N(T_{\rme}) = k, J(T_{\rme}) = j \mid J(0) = i),
\qquad i,j\in\bbM.
\]
It then follows from Lemma~4.1 in \cite{Masu13-ANNOR} that
\begin{equation}
\overline{\vc{A}}(k)\vc{e}
= h \cdot \vc{A}_{\rme} \ast \overline{\vc{D}}(k)\vc{e}
= h \cdot \sum_{l=0}^k \vc{A}_{\rme}(l) \overline{\vc{D}}(k-l)\vc{e}.
\label{eqn-overline{A}(k)e}
\end{equation}
It also follows from Corollary B.1 in \cite{Masu13-ANNOR} that under
the conditions (i) and (ii),
\[
\PP(N(T_{\rme}) > k \mid J(0) = i) \simhm{k} \PP(\lambda T_{\rme} > k).
\]
Thus following the proof of Lemma 3.1 in \cite{Masu09} (see Appendix D
therein), we have for $i,j\in\bbM$,
\begin{equation}
[\overline{\vc{A}}_{\rme}(k)]_{i,j}
= \PP(N(T_{\rme}) > k, J(T_{\rme}) = j \mid J(0) = i)
\simhm{k} [\vc{\pi}]_j \PP(\lambda T_{\rme} >  k).
\label{asymp-overline{A}_e(k)}
\end{equation}
Note here that if $\sqrt{T_{\rme}} \in \calL$, then $\PP(T_{\rme} > x)
= \rme^{-o(\sqrt{x})}$ (see \cite[Lemma~A.1]{Masu13}), where $f(x) =
o(g(x))$ represents $\lim_{x\to\infty}f(x)/g(x)=0$. Thus the
conditions (i) and (ii) imply
\begin{equation}
\lim_{k\to\infty}
{\ooverline{\vc{D}}(k)\vc{e} \over \PP(\lambda T_{\rme} >  k)}
= \vc{0}.
\label{asymp-overline{D}(k)}
\end{equation}
Applying (\ref{asymp-overline{A}_e(k)}), (\ref{asymp-overline{D}(k)}) and
Proposition~\ref{prop-Masu11}~(iii) to (\ref{eqn-overline{A}(k)e}), we
have
\[
\ooverline{\vc{A}}(k)\vc{e}
\simhm{k} h\vc{e}\vc{\pi} 
\sum_{k=0}^{\infty}\overline{\vc{D}}(k)\vc{e} \cdot \PP(\lambda T_{\rme} >  k)
= \rho\vc{e} \cdot \PP(\lambda T_{\rme} >  k),
\]
where the last equality is due to (\ref{defn-lambda}) and
(\ref{defn-rho}). As a result, Assumption~\ref{assumpt-subexp} holds
for $Y = \lambda T_{\rme}$ and $\vc{c}^{\rm A} = \rho\vc{e}$. \qed

\section{Main results}\label{sec-main-result}

In this section, we first provide some basic results on the BMAP/GI/1
retrial queue (subsection~\ref{subsec-model-description}). We then
show a stochastic-decomposition-like result of the stationary queue
length in the BMAP/GI/1 retrial queue
(subsection~\ref{subsec-decomp}). Combining the
stochastic-decomposition-like result with matrix analytic methods, we
prove the subexponential tail equivalence of the stationary queue
length distributions in the BMAP/GI/1 queues with and without retrials
(subsection~\ref{subsec-tail-equiv}).

\subsection{BMAP/GI/1 retrial queue}\label{subsec-model-description}

We begin with the description of the BMAP/GI/1 retrial queue.
Customers arrive at a single-server system with no buffer according to
BMAP $\{\vc{C},\vc{D}(k);k\in\bbN\}$. Such customers are called {\it
  primary customers}. If primary customers arriving in a batch find
the server idle, then one of them immediately occupies the server and
the others join the orbit (virtual waiting line); otherwise all of
them join the orbit. The customers in the orbit are called {\it
  retrial customers}.  We assume that the sojourn times of retrial
customers are i.i.d.\ according to an exponential distribution with
mean $1/\mu$.  A retrial customer tries to occupy the server when its
sojourn time in the orbit expires. If the server is idle, the retrial
customer is accepted; otherwise it goes back to the orbit, i.e.,
becomes a retrial customer again. We also assume that the service
times of primary and retrial customers are independent of BMAP
$\{\vc{C},\vc{D}(k);k\in\bbN\}$ and i.i.d.\ according to a general
distribution function $H$ on $[0,\infty)$ with mean $h \in
  (0,\infty)$.

We now consider the queue length process in the BMAP/GI/1 retrial
queue. As in subsection~\ref{subsec-BMAP/GI/1}, let $J(t)$ ($t \ge 0$)
denote the state of the background Markov chain at time $t$. Let
$Q^{(\mu)}(t)$ ($t \ge 0$) denote the number of retrial customers in
the orbit at time $t$. Further let $S^{(\mu)}(t)$ ($t \ge 0$) denote
the number of customers in the server at time $t$. Clearly,
$S^{(\mu)}(t) \in \{0,1\}$ for all $t \ge 0$ and $\{L^{(\mu)}(t) :=
Q^{(\mu)}(t) + S^{(\mu)}(t);t \ge 0\}$ is the queue length, i.e., the
total number of customers in the server and orbit.

By definition, the process
$\{(S^{(\mu)}(t),Q^{(\mu)}(t),J(t));t\ge0\}$ is a semi-regenerative
process (see \cite[Chapter 10, Section 6]{Cinl75}) such that
regenerative points are service completion instants, i.e., time points
at each of which the service of a customer is completed. Let $0 \le
\tau_0 \le \tau_1 \le \tau_2 \le \cdots $ denote service completion
instants. It then follows that $S^{(\mu)}(\tau_n) = 0$ for all $n \in
\bbZ_+$ and $\{(Q^{(\mu)}(\tau_n),J(\tau_n),\tau_n);n\in\bbZ_+\}$ is a
Markov renewal process (see \cite[Chapter 10, Section 1]{Cinl75}).

\begin{rem}\label{rem-regenerative}
We have a Markov chain by observing
$\{(S^{(\mu)}(t),Q^{(\mu)}(t),J(t));t\ge0\}$ at service beginning
instants, i.e., time points at each of which the service of a customer
starts. Thus service beginning instants can be regenerative points of
$\{(S^{(\mu)}(t),Q^{(\mu)}(t),J(t));t\ge0\}$.
\end{rem}

Recall here that the diagonal elements of $\vc{C}$ are negative and
thus
\begin{eqnarray}
\PP(N(x) = 0, J(x) = i \mid J(0) = i)
= 
[\rme^{\svc{C}x}]_{i,i}>0,
\qquad \forall x > 0,~\forall i \in \bbM. \qquad
\label{add-eqn-02a}
\end{eqnarray}
Note also that $\rme^{(\svc{C} + \svc{D})x} > \vc{O}$ for all $x > 0$
due to the irreducibility of $\vc{C} + \vc{D}$. It then follows from
(\ref{assumpt-lambda>0}) that there exists some $k_0 \in \bbN$ such
that for any $m \in \bbN$,
\begin{eqnarray*}
&&\int_0^{x} \rmd x_m
\int_{0}^{x_m} \rmd x_{m-1} \cdots \int_{0}^{x_2} \rmd x_{1}
\, \rme^{(\svc{C} + \svc{D})x_1} \vc{D}(k_0)
\nonumber
\\
&& \qquad {} \times 
\rme^{(\svc{C} + \svc{D})(x_2-x_1)} \vc{D}(k_0)
 \cdots 
\rme^{(\svc{C} + \svc{D})(x_{m-1} - x_{m-2})} \vc{D}(k_0)
\nonumber
\\
&& \qquad {} \times 
\rme^{(\svc{C} + \svc{D})(x_{m} - x_{m-1})} \vc{D}(k_0)
\rme^{(\svc{C} + \svc{D})(x - x_{m})}
> \vc{O},\qquad \forall x > 0,
\end{eqnarray*}
and thus
\begin{equation}
\PP(N(x) \ge mk_0, J(x) = j \mid J(0) = i) > 0,
\qquad \forall x > 0,~\forall (i,j) \in \bbM^2.
\label{add-eqn-02b}
\end{equation}
It follows from (\ref{add-eqn-02a}) and (\ref{add-eqn-02b}) that the
embedded Markov chain $\{(Q^{(\mu)}(\tau_n),J(\tau_n))\}$ is
irreducible. Further the Markov renewal process
$\{(Q^{(\mu)}(\tau_n),J(\tau_n),\tau_n)\}$ is aperiodic due to the
Markov property of $\{(N(t),J(t))\}$.

It should be noted that for all
$k \in \bbZ_+$ and $i \in \bbM$,
\begin{eqnarray}
d_i(k)
&:=& \EE[\tau_{n+1} \mid Q^{(\mu)}(\tau_n) = k, J(\tau_n)=i]
\nonumber
\\
&\le& [(- \vc{C})^{-1}\vc{D}\vc{e}]_i +  h < \infty. \qquad
\label{mean-renewal-time}
\end{eqnarray}
Therefore, if the embedded Markov chain
$\{(Q^{(\mu)}(\tau_n),J(\tau_n))\}$ is positive recurrent, then for
any initial state, the semi-regenerative process
$\{(S^{(\mu)}(t),Q^{(\mu)}(t),J(t))\}$ has the same limiting
distribution (see \cite[Chapter 10,Theorem~6.12]{Cinl75}). In
addition, if $\rho < 1$, then the embedded Markov chain
$\{(Q^{(\mu)}(\tau_n),J(\tau_n))\}$ is positive recurrent (see, e.g.,
\cite[Theorem~3]{Dudi00}). As a result, if $\rho < 1$, then
$\{(S^{(\mu)}(t),Q^{(\mu)}(t),J(t));t\ge0\}$ is stable (i.e., its
limiting distribution exists; see \cite{Loyn62}) and its limiting
distribution is independent of initial conditions.

On the other hand, if $\rho \ge 1$, then the standard BMAP/GI/1 queue
(without retrials) is unstable \cite{Loyn62}. Thus following the proof
of Theorem 2 in \cite{He00}, we can prove that $\rho < 1$ is a
necessary condition for the stability of the BMAP/GI/1 retrial queue.

The above discussion is summarized in the following:
\begin{lem}\label{lemma-stability}
$\{(S^{(\mu)}(t),Q^{(\mu)}(t),J(t));t\ge0\}$ is stable and its
  limiting distribution is independent of initial conditions if and
  only if $\rho < 1$.
\end{lem}

As stated in subsection~\ref{subsec-BMAP/GI/1}, the stability
condition $\rho < 1$ is assumed.  Thus we define $\vc{p}_0(k)$ and
$\vc{p}_1(k)$ ($k\in\bbZ_+$) as $1 \times M$ vectors such that
\begin{align}
[\vc{p}_0(k)]_i
&= \lim_{t \to \infty}\PP(S^{(\mu)}(t) = 0, Q^{(\mu)}(t)=k, J(t) = i),
& i &\in \bbM,
\label{defn-p_0(k)}
\\
{}
[\vc{p}_1(k)]_i
&= \lim_{t \to \infty}\PP(S^{(\mu)}(t) = 1, Q^{(\mu)}(t)=k, J(t) = i),
& i &\in \bbM,
\label{defn-p_1(k)}
\end{align}
respectively. We also define $\widehat{\vc{p}}_0(z)$ and
$\widehat{\vc{p}}_1(z)$ as
\[
\widehat{\vc{p}}_0(z) = \sum_{k=0}^{\infty}z^k \vc{p}_0(k),
\qquad
\widehat{\vc{p}}_1(z) = \sum_{k=0}^{\infty}z^k \vc{p}_1(k),
\]
respectively.

\begin{rem}\label{rem-p_0(k)}
Since the embedded Markov chain
$\{(Q^{(\mu)}(\tau_n),J(\tau_n));n\in\bbZ_+\}$ is irreducible and
positive recurrent, it has the positive stationary probability
vector. Thus we define $\vc{\varphi}(k)$ ($k \in \bbZ_+$) as a $1
\times M$ vector whose $i$th element $[\vc{\varphi}(k)]_i$
($i\in\bbM$) represents the stationary probability that the embedded
Markov chain $\{(Q^{(\mu)}(\tau_n),J(\tau_n))\}$ is in state $(k,i)$.
It then follows from (\ref{mean-renewal-time}) and Theorem~6.12 in
Chapter 10 of \cite{Cinl75} that
\[
\vc{p}_0(k)
\ge { \vc{\varphi}(k) \int_0^{\infty}\rme^{\svc{C}x} \rmd x
\over \vc{\varphi}\vc{d}}
= { \vc{\varphi}(k) (-\vc{C})^{-1}
\over \vc{\varphi}\vc{d} } > \vc{0},\qquad \forall k \in \bbZ_+,
\]
where $\vc{d}:= (d_i(k))_{(k,i)\in\bbZ_+\times\bbM}$ and
$\vc{\varphi}:=(\vc{\varphi}(0),\vc{\varphi}(1),\vc{\varphi}(2),\dots)$.
Similarly, we can confirm that $\vc{p}_1(k) > \vc{0}$ for all $k \in
\bbZ_+$, though we have to consider another embedded Markov chain of
$\{(S^{(\mu)}(t),Q^{(\mu)}(t),J(t))\}$ observed every time the service
of a customer starts.
\end{rem}

\begin{lem}\label{lem-z-transform-queue-length}
$\widehat{\vc{p}}_0(z)$ and $\widehat{\vc{p}}_1(z)$ satisfy the
  following equations:
\begin{align}
\mu\widehat{\vc{p}}{}_0'(z)\left( z\vc{I} - \widehat{\vc{A}}(z) \right)
&= \widehat{\vc{p}}_0(z)
\left( \vc{C} + z^{-1}\widehat{\vc{D}}(z)\widehat{\vc{A}}(z) \right),
& |z| &< 1,
\label{eqn-hat{p}_0'(z)}
\\
\widehat{\vc{p}}_1(z)\left( z\vc{I} - \widehat{\vc{A}}(z) \right)
&= \widehat{\vc{p}}_0(z)\left( \widehat{\vc{A}}(z) - \vc{I} \right),
& |z| &< 1,
\label{eqn-hat{p}_1(z)}
\end{align}
where $\widehat{\vc{p}}{}_0'(z)=(\rmd /\rmd z)\widehat{\vc{p}}_0(z)$.

\end{lem}

\smallskip

\proof This lemma can be proved in a similar way to that of Theorem~1
in \cite{Kim10b}. However, we here provide a complete proof because
the discussion in Section~\ref{proof-lem-lim-overline{p}_0(k)} uses
some of the symbols introduced to prove this lemma.

We first prove (\ref{eqn-hat{p}_0'(z)}). For this purpose, we consider
a censored process
$\{(\widetilde{Q}^{(\mu)}(t),\widetilde{J}(t));t\ge0\}$ of
$\{(S^{(\mu)}(t),Q^{(\mu)}(t),J(t));t\ge0\}$, which is obtained by
observing $\{(S^{(\mu)}(t),Q^{(\mu)}(t),J(t))\}$ only when
$\{S^{(\mu)}(t) = 0\}$.  It is easy to see that
$\{(\widetilde{Q}^{(\mu)}(t),\widetilde{J}(t))\}$ is a Markov chain
whose transition matrix is given by
\begin{equation}
\widetilde{\vc{T}}
:=
\left(
\begin{array}{ccccc}
\widetilde{\vc{T}}_0(0) & 
\widetilde{\vc{T}}_0(1) &
\widetilde{\vc{T}}_0(2) &
\widetilde{\vc{T}}_0(3) &
\cdots
\\
\widetilde{\vc{T}}_1(-1)& 
\widetilde{\vc{T}}_1(0) &
\widetilde{\vc{T}}_1(1) &
\widetilde{\vc{T}}_1(2) &
\cdots
\\
\vc{O}&
\widetilde{\vc{T}}_2(-1)& 
\widetilde{\vc{T}}_2(0) &
\widetilde{\vc{T}}_2(1) &
\cdots
\\
\vc{O}					&
\vc{O}					&
\widetilde{\vc{T}}_3(-1)& 
\widetilde{\vc{T}}_3(0) &
\cdots
\\
\vdots					&
\vdots					&
\ddots					&
\ddots					&
\ddots
\end{array}
\right),
\label{eqn-tilde{T}-00}
\end{equation}
where
\begin{align}
\widetilde{\vc{T}}_n(-1) &= n\mu\vc{A}(0), 
& n &\in \bbN,
\label{eqn-tilde{T}-01}
\\
\widetilde{\vc{T}}_n(0)  &= n\mu\vc{A}(1) + \vc{D}(1)\vc{A}(0)
+ \vc{C} - n\mu\vc{I}, 
& n &\in \bbZ_+,
\label{eqn-tilde{T}-02}
\\
\widetilde{\vc{T}}_n(k)  &= n\mu\vc{A}(k+1) 
+ \sum_{l=1}^{k+1}\vc{D}(l)\vc{A}(k-l+1), 
& n &\in \bbZ_+,~k \in \bbN.
\label{eqn-tilde{T}-03}
\end{align}

Recall here (see Remark~\ref{rem-p_0(k)}) that the embedded Markov
chain $\{(Q^{(\mu)}(\tau_n),J(\tau_n))\}$ is irreducible and positive
recurrent and thus the censored (continuous-time) Markov chain
$\{(\widetilde{Q}^{(\mu)}(t),\widetilde{J}(t))\}$ is irreducible and
recurrent (see \cite[Chapter 8, Definitions 5.1 and
  5.2]{Brem99}). Therefore $\widetilde{\vc{T}}$ has a unique (up to a
multiplicative factor) positive invariant measure (see \cite[Chapter
  8, Theorem 5.1]{Brem99}). Further since
$\lim_{n\to\infty}|[\widetilde{\vc{T}}_n(0)]_{i,i}| = \infty$, the
unique positive invariant measure is normalized and thus there exists
the unique probability vector $\widetilde{\vc{p}}_0$ such that
$\widetilde{\vc{p}}_0\widetilde{\vc{T}}=\vc{0}$ (see \cite[Chapter 8,
  Theorem 5.1]{Brem99}).

On the other hand, Lemma~\ref{lemma-stability} and
Remark~\ref{rem-p_0(k)} imply that the censored Markov chain
$\{(\widetilde{Q}^{(\mu)}(t),\widetilde{J}(t))\}$ has the unique and
positive limiting distribution (independent of initial conditions). As
a result, $\widetilde{\vc{p}}_0$ is the unique limiting distribution
of the Markov chain $\{(\widetilde{Q}^{(\mu)}(t),\widetilde{J}(t))\}$
(see \cite[Chapter 8, Theorems 5.3 and 6.1]{Brem99}).

We partition $\widetilde{\vc{p}}_0$ as
$(\widetilde{\vc{p}}_0(0),\widetilde{\vc{p}}_0(1),\dots)$, where
\[
[\widetilde{\vc{p}}_0(k)]_i
= \lim_{t \to \infty}\PP(\widetilde{Q}^{(\mu)}(t)=k, \widetilde{J}(t) = i) > 0,
\qquad i \in \bbM.
\]
By definition,
there exists some constant $c > 0$ such that
\begin{equation}
\vc{p}_0(k) = c
\widetilde{\vc{p}}_0(k),
\qquad k \in \bbZ_+.
\label{eqn-tilde{p}_0(k)-01}
\end{equation}
Thus from (\ref{eqn-tilde{p}_0(k)-01}), we have
\begin{equation}
(\vc{p}_0(0),\vc{p}_0(1),\dots) \widetilde{\vc{T}} = \vc{0}.
\label{add-eqn-01}
\end{equation}
It follows from (\ref{eqn-tilde{T}-00})--(\ref{eqn-tilde{T}-03}) and
(\ref{add-eqn-01}) that for $k \in \bbZ_+$,
\begin{eqnarray*}
\vc{0}
&=& \sum_{n=0}^{k+1}\vc{p}_0(n)\widetilde{\vc{T}}_n(k-n)
\nonumber
\\
&=& \sum_{n=0}^{k+1} \vc{p}_0(n) n\mu\vc{A}(k-n+1) +	
\sum_{n=0}^k\vc{p}_0(n) \sum_{l=1}^{k-n+1}\vc{D}(l)\vc{A}(k-n-l+1)
\nonumber
\\
&&{} \qquad + \vc{p}_0(k) ( \vc{C} - k\mu\vc{I} ),
\end{eqnarray*}
where the summation over the empty set is defined as zero.
Multiplying both sides of the above equation by $z^k$ and summing them
for all $k\in\bbZ_+$, we obtain
\begin{eqnarray*}
\vc{0}
&=& \mu\sum_{n=1}^{\infty}n z^{n-1}\vc{p}_0(n) 
\sum_{k=n-1}^{\infty}z^{k-n+1}\vc{A}(k-n+1)
\nonumber
\\
&&{} +  \sum_{n=0}^{\infty} z^n \vc{p}_0(n) 	
 \cdot  z^{-1}\sum_{k=n}^{\infty} z^{k-n+1} 
\sum_{l=1}^{k-n+1}\vc{D}(l)\vc{A}(k-n-l+1)
\nonumber
\\
&&{}
+  \sum_{k=0}^{\infty} z^k \vc{p}_0(k) \vc{C} 
- \mu z\sum_{k=1}^{\infty} k z^{k-1} \vc{p}_0(k)
\nonumber
\\
&=& \mu\widehat{\vc{p}}{}_0'(z)
 \left( \widehat{\vc{A}}(z) - z\vc{I} \right)
 + \widehat{\vc{p}}_0(z)
\left( \vc{C} + z^{-1} \widehat{\vc{D}}(z)\widehat{\vc{A}}(z) \right),
\end{eqnarray*}
which leads to (\ref{eqn-hat{p}_0'(z)}). 

Next we prove (\ref{eqn-hat{p}_1(z)}). Let $r_i(k)$ ($k\in\bbZ_+,
i\in\bbM$) denote the stationary probability that the number of
retrial customers is equal to $k$ and the background state is $i$
immediately after the service of a customer starts, which is
well-defined due to Lemma~\ref{lemma-stability}. Note here that the
time-average number of customers in service is equal to the arrival
rate $\lambda$ of primary customers.  It then follows that for $k \in
\bbZ_+$,
\[
\vc{r}(k) 
:=(r_i(k))_{i\in\bbM}
= {1 \over \lambda} 
\left( \vc{p}_0(k+1) (k+1)\mu + \sum_{l=0}^k\vc{p}_0(l)\vc{D}(k-l+1) \right),
\]
which yields
\begin{equation}
\widehat{\vc{r}}(z)
:= \sum_{k=0}^{\infty}z^k \vc{r}(k) 
= {1 \over \lambda} 
\left( \mu \widehat{\vc{p}}{}'_0(z)
+ z^{-1}\widehat{\vc{p}}_0(z)\widehat{\vc{D}}(z) 
\right).
\label{eqn-hat{r}(z)}
\end{equation}
Recall that $\{(S^{(\mu)}(t),Q^{(\mu)}(t),J(t));t\ge0\}$ is a
semi-regenerative process that is regenerated every time the service
of a customer starts (see Remark~\ref{rem-regenerative}). Since the
mean regenerative cycle is equal to $1/\lambda$, it follows from
Theorem~6.12 in Chapter~10 of \cite{Cinl75} that
\[
\widehat{\vc{p}}_1(z)
= \lambda\widehat{\vc{r}}(z)
\int_0^{\infty} e^{(\svc{C} + \widehat{\svc{D}}(z))x}\overline{H}(x) \rmd x,
\]
where $\overline{H}(x) = 1 - H(x)$ for $x \ge 0$.
Combining this with (\ref{eqn-hat{r}(z)}) yields
\begin{eqnarray}
\widehat{\vc{p}}_1(z)
&=& \left( \mu \widehat{\vc{p}}{}_0'(z) 
+ z^{-1}\widehat{\vc{p}}_0(z)\widehat{\vc{D}}(z) \right)
\cdot 
\int_0^{\infty} e^{(\svc{C} + \widehat{\svc{D}}(z)) x} \overline{H}(x) \rmd x
\nonumber
\\
&=& \left( \mu \widehat{\vc{p}}{}_0'(z) 
+ z^{-1}\widehat{\vc{p}}_0(z)\widehat{\vc{D}}(z) \right)
\cdot 
\left( \vc{C} + \widehat{\vc{D}}(z) \right)^{-1}
\left( \widehat{\vc{A}}(z) - \vc{I} \right), \qquad
\label{eqn-01}
\end{eqnarray}
for all $0 \le |z| < 1$. From (\ref{eqn-01}) and
(\ref{eqn-hat{p}_0'(z)}), we have
\begin{eqnarray*}
\widehat{\vc{p}}_1(z)\left( z\vc{I} - \widehat{\vc{A}}(z) \right)
&=& \left( \mu \widehat{\vc{p}}{}_0'(z) 
+ z^{-1}\widehat{\vc{p}}_0(z)\widehat{\vc{D}}(z) \right)
\\
&& {} \quad 
\times \left( \vc{C} + \widehat{\vc{D}}(z) \right)^{-1}
\left( \widehat{\vc{A}}(z) - \vc{I} \right)
\left( z\vc{I} - \widehat{\vc{A}}(z) \right)
\nonumber
\\
&=& \left( \mu \widehat{\vc{p}}{}_0'(z) 
+ z^{-1}\widehat{\vc{p}}_0(z)\widehat{\vc{D}}(z) \right)
\\
&& {} \quad 
\times \left( z\vc{I} - \widehat{\vc{A}}(z) \right)
\left( \vc{C} + \widehat{\vc{D}}(z) \right)^{-1}
\left( \widehat{\vc{A}}(z) - \vc{I} \right)
\nonumber
\\
&=& \widehat{\vc{p}}_0(z)
\left\{ \left( \vc{C} + z^{-1}\widehat{\vc{D}}(z)\widehat{\vc{A}}(z) \right)
+ z^{-1}\widehat{\vc{D}}(z)\left( z\vc{I} - \widehat{\vc{A}}(z) \right)
\right\}
\nonumber
\\
&& {} \quad 
\times \left( \vc{C} + \widehat{\vc{D}}(z) \right)^{-1}
\left( \widehat{\vc{A}}(z) - \vc{I} \right)
\nonumber
\\
&=& \widehat{\vc{p}}_0(z) \left( \widehat{\vc{A}}(z) - \vc{I} \right),
\end{eqnarray*}
where the second equality holds because $\vc{C} + \widehat{\vc{D}}(z)$
and $\widehat{\vc{A}}(z)$ are commutative. \qed

\subsection{Stochastic-decomposition-like result}\label{subsec-decomp}

Let $\vc{x}^{(\mu)}(k)$ ($k\in\bbZ_+$) denote a $1 \times M$ vector
such that
\[
[\vc{x}^{(\mu)}(k)]_i = \lim_{t \to \infty}\PP(L^{(\mu)}(t) = k, J(t) = i),
\qquad i \in \bbM,
\]
where $L^{(\mu)}(t) = Q^{(\mu)}(t)+S^{(\mu)}(t)$.
Further let $\widehat{\vc{x}}{}^{(\mu)}(z) =
\sum_{k=0}^{\infty}z^k\vc{x}^{(\mu)}(k)$. By definition, we then have
\begin{eqnarray}
\vc{x}{}^{(\mu)}(k)
&=& 
\left\{
\begin{array}{ll}
\vc{p}_0(0), & k = 0,
\\
\vc{p}_0(k) + \vc{p}_1(k-1), & k \in \bbN,
\end{array}
\right.
\label{eqn-x^{(mu)}(k)}
\\
\widehat{\vc{x}}{}^{(\mu)}(z) 
&=& \widehat{\vc{p}}_0(z) + z \widehat{\vc{p}}_1(z).
\label{eqn-widehat{x}^{(mu)}(z)}
\end{eqnarray}

The following lemma is an extension of the stochastic decomposition of
the stationary queue length in the M${}^X$/GI/1 retrial queue (see
Proposition 1 in \cite{Yama12}).
\begin{lem}\label{lem-decomposition}

For $\mu \in (0,\infty)$, 
\begin{align}
&&
\widehat{\vc{x}}{}^{(\mu)}(z)
&=
{\widehat{\vc{p}}_0(z) \over 1 - \rho} \cdot \widehat{\vc{X}}(z),
& |z| &< 1, &&
\label{decomp-01}
\\
&&
\vc{g}\widehat{\vc{X}}(z) &= \widehat{\vc{x}}(z),
& |z| &< 1, &&
\label{decomp-02}
\\
&&
\lim_{z\uparrow1}\widehat{\vc{X}}(z) &= \vc{e}\vc{\pi}, &&&&
\label{decomp-03}
\end{align}
where $\widehat{\vc{X}}(z):=\sum_{k=-\infty}^{\infty}z^k \vc{X}(k)$
($|z| < 1$) is defined as
\begin{equation}
\widehat{\vc{X}}(z) = (1-\rho)(z-1)
\left( z\vc{I} - \widehat{\vc{A}}(z) \right)^{-1}\widehat{\vc{A}}(z).
\label{defn-X(z)}
\end{equation}
\end{lem}

\smallskip

\proof Applying
(\ref{eqn-hat{p}_1(z)}) to (\ref{eqn-widehat{x}^{(mu)}(z)}) yields
\begin{eqnarray*}
\widehat{\vc{x}}{}^{(\mu)}(z)
&=& \widehat{\vc{p}}_0(z)
\left\{\vc{I} + z \left( \widehat{\vc{A}}(z) - \vc{I} \right)
\left( z\vc{I} - \widehat{\vc{A}}(z) \right)^{-1} \right\}
\nonumber
\\
&=& \widehat{\vc{p}}_0(z)(z-1)\widehat{\vc{A}}(z)
 \left( z\vc{I} - \widehat{\vc{A}}(z) \right)^{-1}
\nonumber
\\
&=& \widehat{\vc{p}}_0(z)(z-1)
 \left( z\vc{I} - \widehat{\vc{A}}(z) \right)^{-1}\widehat{\vc{A}}(z).
\end{eqnarray*}
Combining this with (\ref{defn-X(z)}), we have
(\ref{decomp-01}). Further (\ref{decomp-02}) follows from
(\ref{eqn-hat{x}^{(infty)}(z)}) and (\ref{defn-X(z)}).

Finally we prove (\ref{decomp-03}). Let $\sigma_i(z)$'s ($z>0$,
$i=1,2,\dots,M$) denote the eigenvalues of $\widehat{\vc{A}}(z)$ such
that $|\sigma_1(z)| \ge |\sigma_2(z)| \ge \cdots \ge |\sigma_M(z)|$.
Note here that since $\vc{A}=\widehat{\vc{A}}(1)$ is irreducible, so
is $\widehat{\vc{A}}(z)$ for $0 < z < r_A$, where $r_A$ is the
convergence radius of $\sum_{k=0}^{\infty}z^k\vc{A}(k)$. Thus
$\sigma_1(z)$ is the Perron-Frobenius eigenvalue and $\sigma_1(z) >
|\sigma_2(z)| \ge \cdots \ge |\sigma_M(z)|$. In addition, from
$\sigma_1(1)=1$, $\vc{\pi}\vc{A}=\vc{\pi}$ and $\vc{A}\vc{e}=\vc{e}$,
we have
\[
\sigma_1'(1) :=
\left. {\rmd \over \rmd z} \sigma_1(z) \right|_{z=1} 
= \vc{\pi} \sum_{k=1}^{\infty} k\vc{A}(k) \vc{e} = \rho.
\]
Therefore, following the proof of Lemma 3.3 in \cite{Kimu10}, we can
show that
\[
\lim_{z\uparrow1}(z-1)
\left( z\vc{I} - \widehat{\vc{A}}(z) \right)^{-1}
= { \vc{e}\vc{\pi} \over 1 - \rho }.
\]
Applying this to (\ref{defn-X(z)}) yields (\ref{decomp-03}). \qed

\medskip

We conclude this subsection with
some remarks on the coefficient matrices $\vc{X}(k)$ ($k\in\bbZ$) of
the power series expansion of $\widehat{\vc{X}}(z)$ in
(\ref{defn-X(z)}).

Combining (\ref{defn-X(z)}) with (\ref{RG-factorization}), 
we have for $|z| < 1$,
\begin{eqnarray}
\widehat{\vc{X}}(z) 
&=& (1-1/z)
(\vc{I} - \vc{G}/z)^{-1} 
\nonumber
\\
&& {} \times 
(1-\rho)(\vc{I}-\vc{U}(0))^{-1}
\left( \vc{I} - \widehat{\vc{R}}(z) \right) ^{-1}\widehat{\vc{A}}(z)
\nonumber
\\
&=& \widehat{\vc{X}}_1(z) \widehat{\vc{X}}_2(z),
\label{eqn-widehat{X}(z)}
\end{eqnarray}
where $\widehat{\vc{X}}_1(z) :=\sum_{k=0}^{\infty}z^{-k}\vc{X}_1(k)$
and $\widehat{\vc{X}}_2(z) :=\sum_{k=0}^{\infty}z^k\vc{X}_2(k)$ are
given by
\begin{align}
\widehat{\vc{X}}_1(z) 
&= (1-1/z) (\vc{I} - \vc{G}/z)^{-1},
& |z| &< 1,
\label{defn-X_1(k)}
\\
\widehat{\vc{X}}_2(z) 
&=
(1-\rho) (\vc{I} - \vc{U}(0))^{-1}
\left( \vc{I} - \widehat{\vc{R}}(z) \right)^{-1}\widehat{\vc{A}}(z),
& |z| &\le 1.
\label{defn-X_2(k)}
\end{align}
From (\ref{eqn-widehat{X}(z)}) and (\ref{defn-X_1(k)}), we have
\begin{eqnarray}
\vc{X}(k)
&=& \sum_{m=\max(-k,0)}^{\infty} \vc{X}_1(m) \vc{X}_2(k+m),
\qquad k \in \bbZ,
\label{eqn-X(k)}
\\
\vc{X}_1(k)
&=& \left\{
\begin{array}{ll}
\vc{I}, & k=0,
\\
\vc{G}^k - \vc{G}^{k-1}, & k \in \bbN.
\end{array}
\right.
\label{eqn-X_1(k)}
\end{eqnarray}
Substituting (\ref{eqn-X_1(k)}) into (\ref{eqn-X(k)}) yields
\begin{eqnarray}
\vc{X}(k)
&=& \vc{X}_2(k)
+ \sum_{m=1}^{\infty} \left( \vc{G}^m -\vc{G}^{m-1} \right)
\vc{X}_2(k+m),
\qquad k \in \bbZ_+.
\label{eqn-X(k)-02}
\end{eqnarray}
Pre-multiplying both sides of (\ref{eqn-X(k)-02}) by $\vc{g}$ and
using (\ref{decomp-02}), we have
\begin{equation}
\vc{g}\vc{X}(k)  = \vc{g} \vc{X}_2(k) = \vc{x}(k),\qquad k \in \bbZ_+.
\label{eqn-gX}
\end{equation}

Equation (\ref{defn-X_2(k)}) implies that $\vc{X}_2(k) \ge \vc{O}$ for
all $k \in \bbZ_+$. On the other hand, (\ref{eqn-X(k)-02}) shows that
$\vc{X}(k)$ ($k\in\bbZ$) itself may not be nonnegative.  It should be
noted that if background state space $\bbM=\{1\}$, i.e., the BMAP/GI/1
retrial queue is reduced to the M${}^X$/GI/1 retrial queue, then
$\vc{g}=1$ and thus (\ref{eqn-gX}) yields
\[
\vc{X}(k) = \vc{X}_2(k) = \vc{x}(k),\qquad k \in \bbZ_+,
\]
which shows that $\{\vc{X}(k)\}$ and $\{\vc{X}_2(k)\}$ are equivalent
to the stationary queue length distribution in the M${}^X$/GI/1
retrial queue.

\subsection{Subexponential tail equivalence}\label{subsec-tail-equiv}

In this subsection, we present the main theorem. To this end, we
provide three lemmas.
\begin{lem}\label{lem-X(-k)}
$\vc{X}(-k) = O(\gamma^k) \vc{e}\vc{e}^\rmt$ for some $\gamma \in
  (0,1)$, where $f(x) = O(g(x))$ represents
  $\limsup_{x\to\infty}|f(x)/g(x)| < \infty$.
\end{lem}

\proof According to Remark~\ref{rem-G}, we fix $\gamma \in
(0,1)$ such that
\begin{equation}
\vc{G}^k = \vc{e}\vc{g} + O(\gamma^k) \vc{e}\vc{e}^\rmt,
\qquad k \in \bbZ_+.
\label{eqn-G-k}
\end{equation}
From (\ref{eqn-X_1(k)}), we then have
\begin{equation}
\vc{X}_1(k) = O(\gamma^k)\vc{e}\vc{e}^\rmt,
\qquad k \in \bbZ_+.
\label{order-X_1(k)}
\end{equation}
From (\ref{defn-X_2(k)}), we also have
$\sum_{l=0}^{\infty}\gamma^l\vc{X}_2(l) < \infty$.  Therefore
substituting (\ref{order-X_1(k)}) to (\ref{eqn-X(k)}) yields for $k
\in \bbN$,
\begin{eqnarray*}
\gamma^{-k}\vc{X}(-k)
&=&  \gamma^{-k}\sum_{m=k}^{\infty} \vc{X}_1(m) \vc{X}_2(-k+m)
\nonumber
\\
&=&  \vc{e}\vc{e}^\rmt 
\sum_{m=k}^{\infty} O(\gamma^{-k+m})\vc{X}_2(-k+m) < \infty,
\end{eqnarray*}
which completes the proof.
\qed

\begin{lem}\label{lem-lim-overline{X}(k)}
If Assumption~\ref{assumpt-subexp} holds, then
\begin{equation}
\overline{\vc{X}}(k) 
\simhm{k} \vc{e} \overline{\vc{x}}(k)
\simhm{k} { \vc{\pi}\vc{c}^{\rm A} \over 1 - \rho}\vc{e}\vc{\pi}
\cdot \PP(Y>k). 
\label{limit-X_2(k)}
\end{equation}
\end{lem}

\proof From (\ref{eqn-X(k)-02}), we have
\begin{eqnarray}
\overline{\vc{X}}(k)
&=& \overline{\vc{X}}_2(k)
+ \sum_{m=1}^{\infty} \left( \vc{G}^m -\vc{G}^{m-1} \right)
\overline{\vc{X}}_2(k+m)
\nonumber
\\
&=& \sum_{m=0}^{\infty} \vc{G}^m
\left\{ \overline{\vc{X}}_2(k+m) - \overline{\vc{X}}_2(k+m+1) \right\}
\nonumber
\\
&=& \sum_{m=0}^{\infty} \vc{G}^m \vc{X}_2(k+m+1) \ge \vc{O},
\qquad k \in \bbZ_+.
\label{eqn-overline{X}(k)-02}
\end{eqnarray}
Recall here that $\lim_{m\to\infty}\vc{G}^m = \vc{e}\vc{g}$ (see
Remark~\ref{rem-G}). Thus for any $\varepsilon > 0$, there exists some
$m_0 :=m_0(\varepsilon) \in \bbN$ such that for all $m \ge m_0$,
\begin{equation}
(1-\varepsilon)\vc{e}\vc{g}
\le \vc{G}^m
\le (1+\varepsilon)\vc{e}\vc{g}.
\label{ineqn-G^m}
\end{equation}
Applying (\ref{ineqn-G^m}) to
(\ref{eqn-overline{X}(k)-02}) and using (\ref{eqn-gX}), we obtain
\begin{eqnarray}
\liminf_{k\to\infty}
{\overline{\vc{X}}(k) \over \PP(Y>k)}
&\ge& \liminf_{k\to\infty}\sum_{m=m_0}^{\infty}\vc{G}^m
{\vc{X}_2(k+m+1) \over \PP(Y>k)}
\nonumber
\\
&\ge&  (1-\varepsilon)\vc{e}
\cdot \liminf_{k\to\infty}{\overline{\vc{x}}(k+m_0) \over \PP(Y>k)}.
\label{asymp-X_2(k)-01}
\end{eqnarray}
Similarly,
\begin{eqnarray}
\limsup_{k\to\infty}
{\overline{\vc{X}}(k) \over \PP(Y>k)}
&\le& (1+\varepsilon)\vc{e}
\cdot \limsup_{k\to\infty}{\overline{\vc{x}}(k+m_0) \over \PP(Y>k)}
\nonumber
\\
&& {} \qquad
+  \sum_{m=0}^{m_0-1}\vc{G}^m
\limsup_{k\to\infty}{\vc{X}_2(k+m+1) \over \PP(Y>k)}.
\label{asymp-X_2(k)-02a}
\end{eqnarray}
Since $\vc{g} > \vc{0}$, there exists some constant $K > 0$ such that
$\vc{G}^m \le K \vc{e}\vc{g}$ for all $m \in \bbZ_+$. Therefore from
(\ref{asymp-X_2(k)-02a}), we have
\begin{eqnarray}
\limsup_{k\to\infty}
{\overline{\vc{X}}(k) \over \PP(Y>k)}
&\le& (1+\varepsilon)\vc{e}
\cdot \limsup_{k\to\infty}{\overline{\vc{x}}(k+m_0) \over \PP(Y>k)}
\nonumber
\\
&& {} 
+ K\vc{e} \cdot \sum_{m=0}^{m_0-1}
\limsup_{k\to\infty}{\overline{\vc{x}}(k+m) - \overline{\vc{x}}(k+m+1) 
\over \PP(Y>k)}. \qquad
\label{asymp-X_2(k)-02b}
\end{eqnarray}

It follows from Assumption~\ref{assumpt-subexp} and $Y \in
\calS \subset \calL$ that for any fixed $m \in \bbZ_+$,
\begin{eqnarray}
\lim_{k\to\infty}{\overline{\vc{x}}(k+m) \over \PP(Y > k)}
= \lim_{k\to\infty}{\overline{\vc{x}}(k+m+1) \over \PP(Y > k)}
= { \vc{\pi}\vc{c}^{\rm A} \over 1 - \rho}\vc{\pi}.
\label{lim-X_2(k)-02}
\end{eqnarray}
Substituting (\ref{lim-X_2(k)-02}) into (\ref{asymp-X_2(k)-01}) and
(\ref{asymp-X_2(k)-02b}) and letting $\varepsilon \downarrow 0$, we
obtain
\[
\lim_{k\to\infty}{\overline{\vc{X}}(k) \over \PP(Y > k)}
= { \vc{\pi}\vc{c}^{\rm A} \over 1 - \rho}\vc{e}\vc{\pi}.
\]
Combining this with Proposition~\ref{prop-Masu13-ANNOR} yields
(\ref{limit-X_2(k)}). \qed

\begin{lem}\label{lem-lim-overline{p}_0(k)}
If Assumption~\ref{assumpt-subexp} holds, then
$\lim_{k\to\infty}\overline{\vc{p}}_0(k) / \PP(Y > k) = \vc{0}$.
\end{lem}

Lemma~\ref{lem-lim-overline{p}_0(k)} is key to the proof of the main
theorem (Theorem~\ref{thm-main} below). We postpone, however, the
proof of this lemma until the next section because the proof is
somewhat long and technical.

The main theorem of this paper is as follows:
\begin{thm}\label{thm-main}
If Assumption~\ref{assumpt-subexp} holds, then
$\overline{\vc{x}}^{(\mu)}(k) \simhm{k} \overline{\vc{x}}(k)$.
\end{thm}

\proof From (\ref{decomp-01}), we have
\begin{eqnarray}
\lefteqn{
(1-\rho)\overline{\vc{x}}^{(\mu)}(k)
}
~~&&
\nonumber
\\
&=& \widehat{\vc{p}}_0(1)\overline{\vc{X}}(k)
+ \sum_{m=-\infty}^k \overline{\vc{p}}_0(k-m)\vc{X}(m)
\nonumber
\\
&=& \widehat{\vc{p}}_0(1)\overline{\vc{X}}(k)
+ \sum_{m=1}^{\infty} \overline{\vc{p}}_0(k+m)\vc{X}(-m)
+ \sum_{m=0}^k \overline{\vc{p}}_0(k-m)\vc{X}(m). \quad
\label{3terms-overline{X}(k)}
\end{eqnarray}
Lemma~\ref{lem-lim-overline{X}(k)} implies that
\begin{equation}
\lim_{k\to\infty}
\widehat{\vc{p}}_0(1){\overline{\vc{X}}(k) \over \PP(Y > k)}
= \widehat{\vc{p}}_0(1) \vc{e}
\cdot {\vc{\pi}\vc{c}^{\rm A} \over 1 - \rho}
\vc{\pi}.
\label{lim-1st-term}
\end{equation}
It follows from (\ref{decomp-01}), (\ref{decomp-03}) and
$\widehat{\vc{x}}{}^{(\mu)}(1) = 1$ that
\begin{equation}
\widehat{\vc{p}}_0(1) \vc{e} = 1 - \rho.
\label{eqn-widehat{p}_0(1)}
\end{equation}
Substituting (\ref{eqn-widehat{p}_0(1)}) into (\ref{lim-1st-term})
yields
\begin{equation}
\lim_{k\to\infty}
\widehat{\vc{p}}_0(1){\overline{\vc{X}}(k) \over \PP(Y > k)}
= \vc{\pi}\vc{c}^{\rm A}
\cdot 
\vc{\pi}.
\label{lim-1st-term-02}
\end{equation}
Further since $\{\overline{\vc{p}}_0(k);k\in\bbZ_+\}$ is
nonincreasing, Lemmas~\ref{lem-X(-k)} and
\ref{lem-lim-overline{p}_0(k)} imply that
\begin{equation}
\limsup_{k\to\infty}
\left|
\sum_{m=1}^{\infty} {\overline{\vc{p}}_0(k+m) \over \PP(Y>k) }\vc{X}(-m)
\right|
\le
\limsup_{k\to\infty}{\overline{\vc{p}}_0(k) \over \PP(Y>k) }
\sum_{m=1}^{\infty} |\vc{X}(-m)|
= \vc{0}.
\label{lim-2nd-term}
\end{equation}
Applying (\ref{lim-1st-term-02}) and (\ref{lim-2nd-term}) to
(\ref{3terms-overline{X}(k)}), we have
\[
\lim_{k\to\infty}
{\overline{\vc{x}}^{(\mu)}(k) \over \PP(Y > k) }
= {\vc{\pi}\vc{c}^{\rm A} \over 1 - \rho}  \vc{\pi}
+ {1 \over 1 - \rho}
\sum_{m=0}^k { \overline{\vc{p}}_0(k-m) \vc{X}(m) \over \PP(Y > k) }. 
\]
Therefore to complete the proof, it suffices to show that
\[
\limsup_{k\to\infty} 
\left| \sum_{m=0}^k 
{ \overline{\vc{p}}_0(k-m) \vc{X}(m) \over \PP(Y > k) } 
\right|
= \vc{O}.
\]
In what follows, we prove this equation. 

According to (\ref{eqn-G-k}) and $\vc{g} > \vc{0}$, there exist some
$K > 0$ and $\gamma \in (0,1)$ such that
\begin{align*}
&&
| \vc{G}^m -\vc{G}^{m-1}| &\le K \gamma^m \vc{e}\vc{g},
& \forall m &\in \bbN, &&
\\
&&
\vc{X}_2(k) &\le K \vc{e}\vc{g}\vc{X}_2(k),
& \forall k &\in \bbZ_+. &&
\end{align*}
Substituting these inequalities into (\ref{eqn-X(k)-02}) and
using (\ref{eqn-gX}) yield for $k \in \bbZ_+$,
\begin{eqnarray*}
|\vc{X}(k)| 
&\le& K\vc{e}
\left(\vc{g}\vc{X}_2(k) + \sum_{m=1}^{\infty} \gamma^m \vc{g}\vc{X}_2(k+m)
\right)
\nonumber
\\
&=& K\vc{e} 
\left( 
\vc{x}(k) + \sum_{m=1}^{\infty} \gamma^m \vc{x}(k+m) 
\right) =: \vc{X}^+(k).
\end{eqnarray*}
Since $\{\overline{\vc{x}}(k);k\in\bbZ_+\}$ is nonincreasing, it
follows from Proposition~\ref{prop-Masu13-ANNOR} that
\[
\sum_{m=1}^{\infty} \gamma^m { \overline{\vc{x}}(k+m) \over \PP(Y>k) }
\le \sup_{k\in\bbZ_+}{ \overline{\vc{x}}(k) \over \PP(Y>k) }
{\gamma \over 1 - \gamma}
< \infty.
\]
Thus using the dominated convergence theorem,
Proposition~\ref{prop-Masu13-ANNOR} and $Y \in \calS \subset \calL$,
we obtain
\begin{eqnarray}
\lim_{k\to\infty} 
{\overline{\vc{X}}{}^+(k) \over \PP(Y>k) }
&=& K\vc{e} 
\left( 
\lim_{k\to\infty} {\overline{\vc{x}}(k) \over \PP(Y>k)}
+ \sum_{m=1}^{\infty} \gamma^m 
\lim_{k\to\infty} {\overline{\vc{x}}(k+m) \over \PP(Y>k)}
\right)
\nonumber
\\
&=& K 
\left({1 + {\gamma \over 1 - \gamma}} \right)
{\vc{\pi}\vc{c}^{\rm A} \over 1 - \rho } \vc{e} \vc{\pi}
= {K \over 1 - \gamma}
{\vc{\pi}\vc{c}^{\rm A} \over 1 - \rho } \vc{e} \vc{\pi}
< \infty. \qquad
\label{lim-overline{X}^+(k)}
\end{eqnarray}
Combining this with Lemma~\ref{lem-lim-overline{p}_0(k)} and
Proposition~\ref{prop-Masu11}, we have
\begin{equation}
\lim_{k\to\infty} {\overline{\vc{p}_0 \ast \vc{X}^+}(k) \over \PP(Y>k)}
= {K \over 1 - \gamma}
{\vc{\pi}\vc{c}^{\rm A} \over 1 - \rho } \widehat{\vc{p}}_0(1)\vc{e} \vc{\pi}
= {K \over 1 - \gamma}
\vc{\pi}\vc{c}^{\rm A} \cdot \vc{\pi},
\label{lim-overline{p_0*X^+}(k)}
\end{equation}
where we use (\ref{eqn-widehat{p}_0(1)}) in the second equality. Note
here that
\[
{ \overline{\vc{p}_0 \ast \vc{X}^+}(k) \over  \PP(Y>k) }
=
{ \widehat{\vc{p}}_0(1)\overline{\vc{X}}{}^+(k) \over  \PP(Y>k) }
+ 
\sum_{m=0}^k {\overline{\vc{p}}_0(k-m)\vc{X}^+(m) \over \PP(Y>k)},
\]
where the first term converges to the right hand side of
(\ref{lim-overline{p_0*X^+}(k)}) as $k \to \infty$, due to
(\ref{lim-overline{X}^+(k)}). Therefore
\[
\lim_{k\to\infty}
\sum_{m=0}^k {\overline{\vc{p}}_0(k-m)\vc{X}^+(m) \over \PP(Y>k)} = \vc{0},
\]
which leads to
\begin{eqnarray*}
\lefteqn{
\limsup_{k\to\infty}
\left|
\sum_{m=0}^k {\overline{\vc{p}}_0(k-m)\vc{X}(m) \over \PP(Y>k)}
\right|
}
\qquad &&
\nonumber
\\
&\le&  \limsup_{k\to\infty}
\sum_{m=0}^k {\overline{\vc{p}}_0(k-m)\vc{X}^+(m) \over \PP(Y>k)}
= \vc{0}.
\end{eqnarray*}
\qed

\medskip

As mentioned in the introduction, Yamamuro~\cite{Yama12} proves the
subexponential tail equivalence of the queue length distributions in
the M${}^X$/GI/1 retrial queues with and without retrials, under the
assumption that the batch size distribution is light-tailed. On the
other hand, Proposition~\ref{lem-assumpt-subexp} shows that
Assumption~\ref{assumpt-subexp} and thus Theorem~\ref{thm-main} do not
necessarily require that $\{\vc{D}(k)\}$ is light-tailed.

\section{Proof of a Key Lemma (Lemma~\ref{lem-lim-overline{p}_0(k)})}
\label{proof-lem-lim-overline{p}_0(k)}

To facilitate the discussion, we apply a change of measure to $\vc{p}_0 := (\vc{p}_0(0),\vc{p}_0(1),\vc{p}_0(2),\dots)$. 
Let $\vc{q}:= (\vc{q}(0),\vc{q}(1),\vc{q}(2),\dots)$ denote a
probability vector
\begin{equation}
\vc{q}(k)
= {\max(k,1)\vc{p}_0(k) 
\over \sum_{l=0}^{\infty}\max(l,1)\vc{p}_0(l)\vc{e}} > \vc{0},
\qquad k \in \bbZ_+,
\label{defn-q(k)}
\end{equation}
where the positivity of $\vc{q}(k)$ follows from Remark~\ref{rem-p_0(k)}.
We then have
\begin{equation}
\vc{q} = {\vc{p}_0\vc{\Delta}^{-1} 
\over \vc{p}_0\vc{\Delta}^{-1}\vc{e}},
\label{eqn-q}
\end{equation}
where $\vc{\Delta}$ is a diagonal matrix such that
\[
\vc{\Delta}
=
\left(
\begin{array}{cccccc}
\vc{I} & 
\vc{O} &
\vc{O} &
\vc{O} &
\vc{O} &
\cdots
\\
\vc{O} & 
\vc{I} &
\vc{O} &
\vc{O} &
\vc{O} &
\cdots
\\
\vc{O} &
\vc{O} & 
\displaystyle{1 \over 2}\vc{I} &
\vc{O} &
\vc{O} &
\cdots
\\
\vc{O} &
\vc{O} &
\vc{O} & 
\displaystyle{1 \over 3}\vc{I} &
\vc{O} &
\cdots
\\
\vc{O} &
\vc{O} &
\vc{O} & 
\vc{O} &
\displaystyle{1 \over 4}\vc{I} &
\ddots
\\
\vdots &
\vdots &
\vdots &
\vdots &
\ddots &
\ddots
\end{array}
\right).
\]

It follows from (\ref{eqn-q}) and $\vc{p}_0\widetilde{\vc{T}}=\vc{0}$
(see (\ref{add-eqn-01})) that $\vc{q} > \vc{0}$ is the stationary
probability vector of the following infinitesimal generator:
\begin{equation}
\breve{\vc{T}} :=
\vc{\Delta}\widetilde{\vc{T}}
=
\left(
\begin{array}{ccccc}
\breve{\vc{T}}_0(0) & 
\breve{\vc{T}}_0(1) &
\breve{\vc{T}}_0(2) &
\breve{\vc{T}}_0(3) &
\cdots
\\
\breve{\vc{T}}_1(-1)& 
\breve{\vc{T}}_1(0) &
\breve{\vc{T}}_1(1) &
\breve{\vc{T}}_1(2) &
\cdots
\\
\vc{O}&
\breve{\vc{T}}_2(-1)& 
\breve{\vc{T}}_2(0) &
\breve{\vc{T}}_2(1) &
\cdots
\\
\vc{O}					&
\vc{O}					&
\breve{\vc{T}}_3(-1)& 
\breve{\vc{T}}_3(0) &
\cdots
\\
\vdots					&
\vdots					&
\ddots					&
\ddots					&
\ddots
\end{array}
\right),
\label{breve-T}
\end{equation}
where
\begin{align}
\breve{\vc{T}}_0(0) 
&= \vc{C} + \vc{D}(1)\vc{A}(0), 
&&
\label{eqn-breve{T}-00a}
\\
\breve{\vc{T}}_0(k)  
&= \sum_{l=1}^{k+1}\vc{D}(l)\vc{A}(k-l+1), 
& k &\in \bbN,
\label{eqn-breve{T}-00b}
\end{align}
and for $n \in \bbN$,
\begin{align}
\breve{\vc{T}}_n(-1) 
&= \mu\vc{A}(0), 
&&
\label{eqn-breve{T}-01}
\\
\breve{\vc{T}}_n(0)  
&= - \mu\vc{I} + \mu\vc{A}(1)
+ {1 \over n}\{ \vc{C} + \vc{D}(1)\vc{A}(0) \}, 
&&
\label{eqn-breve{T}-02}
\\
\breve{\vc{T}}_n(k)  
&=   \mu\vc{A}(k+1) 
+ {1 \over n} \sum_{l=1}^{k+1}\vc{D}(l)\vc{A}(k-l+1), 
& k &\in \bbN.
\label{eqn-breve{T}-03}
\end{align}

For convenience, we uniformize the transition rate matrix
$\breve{\vc{T}}$ as follows:
\begin{equation}
\breve{\vc{P}}:=\vc{I} + {\breve{\vc{T}} \over \mu + \theta},
\label{defn-breve{P}}
\end{equation}
where $\theta = \max_{i\in\bbM}|[\vc{C}]_{i,i}|$.  From
(\ref{breve-T})--(\ref{eqn-breve{T}-03}), we have
\begin{equation}
\breve{\vc{P}}
=
\left(
\begin{array}{ccccc}
\breve{\vc{A}}_0(0) & 
\breve{\vc{A}}_0(1) &
\breve{\vc{A}}_0(2) &
\breve{\vc{A}}_0(3) &
\cdots
\\
\breve{\vc{A}}_1(-1)& 
\breve{\vc{A}}_1(0) &
\breve{\vc{A}}_1(1) &
\breve{\vc{A}}_1(2) &
\cdots
\\
\vc{O}&
\breve{\vc{A}}_2(-1)& 
\breve{\vc{A}}_2(0) &
\breve{\vc{A}}_2(1) &
\cdots
\\
\vc{O}					&
\vc{O}					&
\breve{\vc{A}}_3(-1)& 
\breve{\vc{A}}_3(0) &
\cdots
\\
\vdots					&
\vdots					&
\ddots					&
\ddots					&
\ddots
\end{array}
\right),
\label{eqn-breve{P}}
\end{equation}
where
\begin{align}
\breve{\vc{A}}_0(0) 
&= \vc{I} + {1 \over \mu+\theta} 
\{ \vc{C} + \vc{D}(1)\vc{A}(0) \}, 
&&
\nonumber
\\
\breve{\vc{A}}_0(k)  
&= {1 \over \mu+\theta}\sum_{l=1}^{k+1}\vc{D}(l)\vc{A}(k-l+1), 
& k &\in \bbN,
\label{eqn-breve{A}-00b}
\end{align}
and for $n \in \bbN$,
\begin{align}
\breve{\vc{A}}_n(-1) 
&= {\mu \over \mu+\theta}\vc{A}(0), 
&&
\label{eqn-breve{A}-01}
\\
\breve{\vc{A}}_n(0)  
&= {1 \over \mu+\theta}
\left[
\theta\vc{I} + \mu\vc{A}(1) 
+ {1 \over n} \{ \vc{C} + \vc{D}(1)\vc{A}(0) \}
\right], 
\label{eqn-breve{A}-02}
\\
\breve{\vc{A}}_n(k)  
&= {1 \over \mu+\theta}
\left[
\mu\vc{A}(k+1) 
+ {1 \over n} \sum_{l=1}^{k+1}\vc{D}(l)\vc{A}(k-l+1) 
\right], 
&  k &\in \bbN.
\label{eqn-breve{A}-03}
\end{align}
Let $\breve{\vc{A}}(k)$ ($k\ge -1$) denote
\begin{eqnarray}
\breve{\vc{A}}(-1) 
&=& \breve{\vc{A}}_n(-1)  = {\mu \over \mu + \theta}\vc{A}(0),
\label{eqn-breve{A}(-1)}
\\
\breve{\vc{A}}(0) 
&=& {\theta \over \mu + \theta}\vc{I} + {\mu \over \mu + \theta}\vc{A}(1),
\label{eqn-breve{A}(0)}
\\
\breve{\vc{A}}(k) &=& {\mu \over \mu + \theta}\vc{A}(k+1), \qquad k \in \bbN.
\label{eqn-breve{A}(k)}
\end{eqnarray}
It then follows from (\ref{eqn-breve{A}-01})--(\ref{eqn-breve{A}(k)})
that
\begin{eqnarray}
\lim_{n\to\infty}
\sum_{k=-1}^{\infty}
|\breve{\vc{A}}_n(k) - \breve{\vc{A}}(k)|
\le {1 \over \mu + \theta} \lim_{n\to\infty}{1 \over n} 
\left(|\vc{C}| + \vc{D}\vc{A} \right)
= \vc{O},
\label{lim-breve{A}_n(k)-02}
\end{eqnarray}
and thus
\begin{equation}
\lim_{n\to\infty} \breve{\vc{A}}_n(k) 
= \breve{\vc{A}}(k), \qquad \mbox{uniformly all $k \ge -1$.}
\label{lim-breve{A}_n(k)-01}
\end{equation}

By definition, $\vc{q}= (\vc{q}(0),\vc{q}(1),\vc{q}(2),\dots) >
\vc{0}$ is the stationary probability vector of $\breve{\vc{P}}$. 
Note here that (\ref{breve-T}) and (\ref{defn-breve{P}}) yield
\[
\breve{\vc{P}} 
= \vc{I} + {\vc{\Delta}\widetilde{\vc{T}} \over \mu + \theta}.
\]
Note also that $\widetilde{\vc{T}}$ has the unique stationary
probability vector $\widetilde{\vc{p}}_0 = c^{-1}\vc{p}_0$ (this
equality is due to (\ref{eqn-tilde{p}_0(k)-01})). Therefore $\vc{q}$
is the unique and positive stationary probability vector of
$\breve{\vc{P}}$, which implies that $\breve{\vc{P}}$ is irreducible
and positive recurrent (see \cite[Chapter 3, Theorem 3.1]{Brem99}).
 
As shown in (\ref{eqn-breve{P}}) and (\ref{lim-breve{A}_n(k)-01}),
$\breve{\vc{P}}$ is a level-dependent M/G/1-type stochastic matrix
with asymptotic level-independence. Utilizing this special structure
of $\breve{\vc{P}}$, we can prove Lemma~\ref{lem-asymp-q(k)} below.
\begin{lem}\label{lem-asymp-q(k)}
Suppose Assumption~\ref{assumpt-subexp} holds. The following then
hold.
\begin{enumerate}
\item $\limsup_{k\to\infty}\overline{\vc{q}}(k) / \PP(Y > k)$ is
  finite; and
\item $\lim_{k\to\infty}\overline{\vc{q}}(k) / \PP(Y > k)= c\vc{\pi}$
  for some $c \in (0,\infty)$ if
  $\limsup_{k\to\infty}\ooverline{\vc{D}}(k)\vc{e} / \PP(Y > k)$
  exists.
\end{enumerate}

\end{lem}

From (\ref{defn-q(k)}), we have
\begin{equation}
\vc{p}_0(k) = \vc{p}_0\vc{\Delta}^{-1}\vc{e}\cdot
{\vc{q}(k) \over k}, \qquad k \in \bbN,
\label{eqn-p_0(k)-q(k)}
\end{equation}
which leads to
\[
\overline{\vc{p}}_0(k) 
= \vc{p}_0\vc{\Delta}^{-1}\vc{e}
\cdot \sum_{l=k+1}^{\infty}{\vc{q}(l) \over l}
\le \vc{p}_0\vc{\Delta}^{-1}\vc{e} 
\cdot {\overline{\vc{q}}(k) \over k}, 
\qquad k \in \bbN.
\]
Thus Lemma~\ref{lem-lim-overline{p}_0(k)} is immediate
from statement (i) of Lemma~\ref{lem-asymp-q(k)}. 

\begin{rem}
Although statement (ii) of Lemma~\ref{lem-asymp-q(k)} is not necessary
for Lemma~\ref{lem-lim-overline{p}_0(k)}, the statement is, as far as
we know, the first reported result on the subexponential asymptotics
of level-dependent structured Markov chain. In addition, from
statement (ii), we can guess that under
Assumption~\ref{assumpt-subexp} and additional conditions, the
following locally subexponential asymptotic formula holds:
\begin{equation}
\vc{q}(k) 
\simhm{k} c\vc{\pi} \cdot \PP(Y = k).
\label{locally-asymp-q(k)}
\end{equation}
This can be proved by extending the results on the locally
subexponential asymptotics of {\it level-independent} GI/G/1-type
Markov chains (see Section 4 in \cite{Kimu13}) to a {\it
  level-dependent} M/G/1-type Markov chain with asymptotic
level-independence. If (\ref{locally-asymp-q(k)}) holds, then
(\ref{eqn-p_0(k)-q(k)}) yields
\[
\vc{p}_0(k) 
\simhm{k} c'\vc{\pi} \cdot k^{-1}\PP(Y = k)
\quad \mbox{for some constant $c'>0$}.
\]
Further if the stronger conditions in
Proposition~\ref{lem-assumpt-subexp} are assumed instead of
Assumption~\ref{assumpt-subexp} (of course, other additional
conditions are needed for the locally subexponential asymptotics),
then
\[
\vc{p}_0(k) 
\simhm{k} c'h\vc{\pi} \cdot k^{-1}\PP(\lambda T > k),
\]
where $T$ is the service time.
\end{rem}

We need several technical lemmas to prove
Lemma~\ref{lem-asymp-q(k)}. In the rest of this section, we present
the technical lemmas and then give the proof of
Lemma~\ref{lem-asymp-q(k)} at the end of this section.

We begin with the following lemma.
\begin{lem}\label{overline{breve{A}}(k)>O}
(i) The diagonal elements of $\vc{A}(0)$ and $\breve{\vc{A}}(-1)$ are
  positive; and (ii) for all $k \in \bbZ_+$,
  $\sum_{l=k}^{\infty}\vc{A}(l) > \vc{O}$ and
  $\sum_{l=k-1}^{\infty}\breve{\vc{A}}(l) > \vc{O}$.
\end{lem}

\proof It suffices to prove statements (i) and (ii) for
$\{\vc{A}(k)\}$ due to
(\ref{eqn-breve{A}(-1)})--(\ref{eqn-breve{A}(k)}). It follows from
(\ref{add-eqn-02a}) and (\ref{add-eqn-02b}) that
\[
[\vc{A}(0)]_{i,i} 
= \left[ 
\int_0^{\infty} \rme^{\svc{C} x} \rmd H(x) 
\right]_{i,i} > 0,
\qquad \forall i \in \bbM,
\]
and that there exists some $k_0 \in \bbN$ such that for all $m \in
\bbN$ and $(i,j) \in \bbM^2$,
\begin{eqnarray*}
\left[
\sum_{l=mk_0}^{\infty}\vc{A}(l)
\right]_{i,j}
= \int_0^{\infty} \rmd H(x) \PP(N(x) \ge mk_0, J(x) = j \mid J(0) = i)
> 0,
\end{eqnarray*}
which completes the proof. \qed

For further discussion, we introduce some symbols.
We define $\breve{\vc{P}}_n$ ($n\in\bbN$) as a submatrix of
$\breve{\vc{P}}$ in (\ref{eqn-breve{P}}) such that
\begin{equation}
\breve{\vc{P}}_n
=
\left(
\begin{array}{ccccc}
\breve{\vc{A}}_n(0) & 
\breve{\vc{A}}_n(1) &
\breve{\vc{A}}_n(2) &
\breve{\vc{A}}_n(3) &
\cdots
\\
\breve{\vc{A}}_{n+1}(-1)& 
\breve{\vc{A}}_{n+1}(0) &
\breve{\vc{A}}_{n+1}(1) &
\breve{\vc{A}}_{n+1}(2) &
\cdots
\\
\vc{O}&
\breve{\vc{A}}_{n+2}(-1)& 
\breve{\vc{A}}_{n+2}(0) &
\breve{\vc{A}}_{n+2}(1) &
\cdots
\\
\vc{O}					&
\vc{O}					&
\breve{\vc{A}}_{n+3}(-1)& 
\breve{\vc{A}}_{n+3}(0) &
\cdots
\\
\vdots					&
\vdots					&
\ddots					&
\ddots					&
\ddots
\end{array}
\right).
\label{eqn-breve{P}_n}
\end{equation}
It follows from (\ref{lim-breve{A}_n(k)-02}) and
(\ref{eqn-breve{P}_n}) that
\begin{equation}
\lim_{n\to\infty}\breve{\vc{P}}_n
=
\left(
\begin{array}{ccccc}
\breve{\vc{A}}(0) & 
\breve{\vc{A}}(1) &
\breve{\vc{A}}(2) &
\breve{\vc{A}}(3) &
\cdots
\\
\breve{\vc{A}}(-1)& 
\breve{\vc{A}}(0) &
\breve{\vc{A}}(1) &
\breve{\vc{A}}(2) &
\cdots
\\
\vc{O}&
\breve{\vc{A}}(-1)& 
\breve{\vc{A}}(0) &
\breve{\vc{A}}(1) &
\cdots
\\
\vc{O}					&
\vc{O}					&
\breve{\vc{A}}(-1)& 
\breve{\vc{A}}(0) &
\cdots
\\
\vdots					&
\vdots					&
\ddots					&
\ddots					&
\ddots
\end{array}
\right)
=: \breve{\vc{P}}_{\infty},
\label{eqn-P^{(infty)}}
\end{equation}
where the convergence is uniform over all the elements.  

Recall that $\breve{\vc{P}}$ is irreducible and positive recurrent and
thus the set of states $\{(m,i);m \ge n,\,i\in\bbM\}$ is not closed
for any $n \in \bbN$. Therefore for any $n \in \bbN$, there exists the
minimal nonnegative inverse of $\vc{I} - \breve{\vc{P}}_n$ (see, e.g.,
\cite[Corollary 2 of Lemma 5.4]{Sene06}), which is denoted by $(\vc{I}
- \breve{\vc{P}}_n)^{-1}$ and given by
\begin{equation*}
\left( \vc{I} - \breve{\vc{P}}_n \right)^{-1} 
= \sum_{m=0}^{\infty} (\breve{\vc{P}}_n)^m.
\end{equation*}

Using the inverse $(\vc{I} - \breve{\vc{P}}_n)^{-1}$, we define some
matrices, which play a role in matrix analytic methods.  Let
$\breve{\vc{N}}_n(0,0)$ denote the $M \times M$ northwest corner of
$(\vc{I} - \breve{\vc{P}}_n)^{-1}$, i.e.,
\begin{equation}
\breve{\vc{N}}_n(0,0) 
= \sum_{m=0}^{\infty}\breve{\vc{P}}{}_n^{(m)}(0,0),
\label{defn-N_n(0,0)}
\end{equation}
where $\breve{\vc{P}}{}_n^{(m)}(0,0)$ is the $M \times M$ northwest
corner of $(\breve{\vc{P}}_n)^m$. Further let $\breve{\vc{G}} _n$ ($n
\in \bbN$) and $\breve{\vc{U}}_n(0)$ ($n \in \bbN$) and $\breve
    {\vc{R}}_n(k)$ ($n \in \bbZ_+$, $k\in\bbN$) denote
\begin{eqnarray}
\breve{\vc{G}}_n
&=& \breve{\vc{N}}_n(0,0)\breve{\vc{A}}_n(-1)
= \sum_{m=0}^{\infty}\breve{\vc{P}}{}_n^{(m)}(0,0)\breve{\vc{A}}_n(-1),
\label{eqn-G_n}
\\
\breve{\vc{U}}_n(0) 
&=& \sum_{k=0}^{\infty} \breve{\vc{A}}_n(k)
\prod_{l=n+k}^{n+1}\breve{\vc{G}}_l,
\label{defn-U_n(0)}
\\
\breve{\vc{R}}_n(k)
&=& \sum_{m=0}^{\infty}\breve{\vc{A}}_n(k+m)
\left(
\prod_{l=n+k+m}^{n+k+1}\breve{\vc{G}}_l
\right)
\breve{\vc{N}}_{n+k}(0,0),
\label{defn-R_n(k)}
\end{eqnarray}
respectively, where for $\nu,\eta \in \bbN$,
\[
\prod_{l=\nu}^{\eta}\breve{\vc{G}}_l 
= \left\{
\begin{array}{ll}
\vc{I}, & \nu < \eta,
\\
\breve{\vc{G}}_{\nu}\breve{\vc{G}}_{\nu-1} 
\cdots \breve{\vc{G}}_{\eta}, & \nu \ge \eta.
\end{array}
\right.
\]

In order to interpret the matrices $\breve{\vc{N}}_n(0,0)$,
$\breve{\vc{G}}_n$,$\breve{\vc{U}}_n(0)$ and $\breve{\vc{R}}_n(k)$, we
consider a discrete-time Markov chain
$\{(\breve{L}_m,\breve{J}_m);m\in\bbZ_+\}$ with state space $\bbZ_+
\times \bbM$ and transition matrix $\breve{\vc{P}}$. For simplicity,
we also define $\bbL(n)$ ($n \in \bbZ_+$) as the set of states
$\{(n,i);i\in\bbM\}$. In this setting, the interpretation of
$\breve{\vc{N}}_n(0,0)$, $\breve{\vc{G}}_n$, $\breve{\vc{U}}_n(0)$ and
$\breve{\vc{R}}_n(k)$ is as follows (see \cite{Zhao98}):
\begin{enumerate}
\item $[\breve{\vc{N}}_n(0,0)]_{i,j}$ represents the conditional
  expected number of visits to state $(n,j)$ before entering
  $\bigcup_{\nu=0}^{n-1}\bbL(\nu)$ given that
  $\{(\breve{L}_m,\breve{J}_m)\}$ starts with state $(n,i)$.
\item $[\breve{\vc{G}}_n]_{i,j}$ represents the conditional
  probability that the first passage time to $\bbL(n-1)$ ends with
  state $(n-1,j)$ given that $\{(\breve{L}_m,\breve{J}_m)\}$ starts
  with state $(n,i)$. Note that during the first passage time to
  $\bbL(n-1)$ from $\bbL_n$, $\{(\breve{L}_m,\breve{J}_m)\}$ does not
  visit any state in $\bigcup_{\nu=0}^{n-2}\bbL(\nu)$ because it is
  skip-free to the left (see, e.g., \cite[Chapter 13]{Lato99}).
\item $[\breve{\vc{U}}_n(0)]_{i,j}$ represents the conditional
  probability that the first passage time to
  $\bigcup_{\nu=0}^n\bbL(\nu)$ ends with state $(n,j)$ given that
  $\{(\breve{L}_m,\breve{J}_m)\}$ starts with state $(n,i)$.
\item $[\breve{\vc{R}}_n(k)]_{i,j}$ represents the conditional
  expected number of visits to state $(n+k,j)$ before entering
  $\bigcup_{\nu=0}^{n+k-1}\bbL(\nu)$ given that
  $\{(\breve{L}_m,\breve{J}_m)\}$ starts with state $(n,i)$.
\end{enumerate}

\begin{lem}\label{lem-RG-bounds}
\hfill
\begin{enumerate}
\item For all $n \in \bbN$, $\breve{\vc{G}}_n$ is stochastic matrix;
  and
\item there exists some $\xi \in (0,1)$ such that
\begin{eqnarray}
\sup_{n\in\bbN}\breve{\vc{U}}_n(0) \vc{e} &\le& \xi \vc{e}, 
~~
\sup_{n\in\bbN}\breve{\vc{N}}_n(0,0)\vc{e} 
\le {1 \over 1 - \xi}\vc{e},
\nonumber
\\
\sup_{n\in\bbZ_+} \sum_{k=1}^{\infty}\breve{\vc{R}}_n(k)\vc{e}
&\le& {1 \over 1 - \xi} {1 \over \mu+\theta}
\left[
\mu\widehat{\vc{A}}{}'(1)\vc{e}
+ \vc{D}\widehat{\vc{A}}{}'(1)\vc{e}
+ \widehat{\vc{D}}{}'(1)\vc{e}
\right], \qquad
\label{ineqn-R_n(k)e}
\end{eqnarray}
\end{enumerate}
where $\widehat{\vc{A}}{}'(z)= (\rmd /\rmd z)\widehat{\vc{A}}(z)$ and
$\widehat{\vc{D}}{}'(z)= (\rmd /\rmd z)\widehat{\vc{D}}(z)$.
\end{lem}

\proof Note that $\breve{\vc{P}}$ and thus
$\{(\breve{L}_m,\breve{J}_m)\}$ are irreducible and positive
recurrent.  Note also that $\{(\breve{L}_m,\breve{J}_m)\}$ is
skip-free to the left. Therefore the probabilistic interpretation of
$\breve{\vc{G}}_n$ implies that statement (i) is true.

Next we prove statement (ii). From (\ref{defn-U_n(0)}),
(\ref{eqn-breve{A}(-1)}) and $\breve{\vc{G}}_n\vc{e} = \vc{e}$, we
obtain
\begin{equation}
\breve{\vc{U}}_n(0) \vc{e}
= \sum_{k=0}^{\infty} \breve{\vc{A}}_n(k)\vc{e}
= \vc{e} - \breve{\vc{A}}(-1)\vc{e}
< \vc{e},\qquad \forall n \in \bbN,
\label{add-ineqn-U_n(0)e}
\end{equation}
where the last inequality holds due to statement (i) of
Lemma~\ref{overline{breve{A}}(k)>O}. According to
(\ref{add-ineqn-U_n(0)e}), there exists some $\xi \in (0,1)$ such that
\begin{equation}
\breve{\vc{U}}_n(0) \vc{e} \le \xi \vc{e},
\qquad \forall n \in \bbN.
\label{ineqn-U_n(0)e}
\end{equation}
In addition, the interpretation of $\breve{\vc{N}}_n(0,0)$ and
$\breve{\vc{U}}_n(0)$ implies that
\begin{equation}
\breve{\vc{N}}_n(0,0) 
= \sum_{m=0}^{\infty} \left( \breve{\vc{U}}_n(0) \right)^m
= \left( \vc{I} - \breve{\vc{U}}_n(0) \right)^{-1}.
\label{eqn-breve{N}_n(0,0)}
\end{equation}
Substituting (\ref{ineqn-U_n(0)e}) into (\ref{eqn-breve{N}_n(0,0)})
yields
\begin{equation}
\breve{\vc{N}}_n(0,0)\vc{e} 
=  \left( \vc{I} - \breve{\vc{U}}_n(0) \right)^{-1}\vc{e}
\le {1 \over 1 - \xi}\vc{e},
\qquad \forall n \in \bbN.
\label{ineqn-N_n(0,0)e}
\end{equation}

It remains to prove (\ref{ineqn-R_n(k)e}).
From (\ref{defn-R_n(k)})
and (\ref{ineqn-N_n(0,0)e}), we have for $n\in\bbZ_+$,
\begin{eqnarray*}
\sum_{k=1}^{\infty}
\breve{\vc{R}}_n(k)\vc{e}
&\le& {1 \over 1 - \xi} \sum_{k=1}^{\infty}
\sum_{m=0}^{\infty}\breve{\vc{A}}_n(k+m)\vc{e}
= {1 \over 1 - \xi} 
\sum_{m=0}^{\infty}\overline{\breve{\vc{A}}}_n(m) \vc{e}.
\end{eqnarray*}
From (\ref{eqn-breve{A}-00b}) and (\ref{eqn-breve{A}-03}), we also
have for $n\in\bbZ_+$,
\begin{eqnarray*}
\sum_{m=0}^{\infty}\overline{\breve{\vc{A}}}_n(m) \vc{e}
&\le& 
\sum_{m=0}^{\infty}\overline{\breve{\vc{A}}}_1(m) \vc{e}
= \sum_{k=0}^{\infty}k\breve{\vc{A}}_1(k) \vc{e}
\le \sum_{k=0}^{\infty}(k+1)\breve{\vc{A}}_1(k) \vc{e}
\nonumber
\\
&=& {1 \over \mu+\theta}
\sum_{k=1}^{\infty}k
\left[
\mu\vc{A}(k) 
+ \sum_{l=1}^{k}\vc{D}(l)\vc{A}(k-l) 
\right]\vc{e}
\nonumber
\\
&=& {1 \over \mu+\theta}
\left[
\mu\widehat{\vc{A}}{}'(1)\vc{e}
+ \vc{D}\widehat{\vc{A}}{}'(1)\vc{e}
+ \widehat{\vc{D}}{}'(1)\vc{e}
\right],
\end{eqnarray*}
which is finite due to (\ref{defn-lambda}), (\ref{eqn-rho}) and $\rho
=\lambda h < 1$. As a result, (\ref{ineqn-R_n(k)e}) holds.  \qed

Using Lemma~\ref{lem-RG-bounds} and the dominated convergence theorem,
we obtain the following:
\begin{lem}\label{lem-limit-G_n}
Let
$\breve{\vc{P}}{}_{\infty}^{(m)}(0,0)=\lim_{n\to\infty}\breve{\vc{P}}{}_n^{(m)}(0,0)$
for $m\in\bbZ_+$, which is the $M \times M$ northwest corner of
$(\breve{\vc{P}}_{\infty})^m$. We then have
\begin{eqnarray}
\lim_{n\to\infty}\breve{\vc{G}}_n
&=& \sum_{m=0}^{\infty}\breve{\vc{P}}{}_{\infty}^{(m)}(0,0)\breve{\vc{A}}(-1)
=: \breve{\vc{G}} > \vc{O},
\label{lim-G_n}
\\
\lim_{n\to\infty}
\breve{\vc{U}}_n(0) 
&=& \sum_{k=0}^{\infty} \breve{\vc{A}}(k)\breve{\vc{G}}{}^k 
=: \breve{\vc{U}}(0)  > \vc{O},
\label{eqn-U^{(infty)}(0)}
\\
\lim_{n\to\infty}
\breve{\vc{N}}_n(0,0) 
&=& \left( \vc{I} - \breve{\vc{U}}(0) \right)^{-1} > \vc{O},
\label{eqn-N^{(infty)}(0)}
\end{eqnarray}
and for $k \in \bbN$,
\begin{eqnarray}
\lim_{n\to\infty}
\breve{\vc{R}}_n(k)
&=& \sum_{m=0}^{\infty}\breve{\vc{A}}(k+m)\breve{\vc{G}}{}^m
\left( \vc{I} - \breve{\vc{U}}(0) \right)^{-1} 
=: \breve{\vc{R}}(k) > \vc{O}.
\label{eqn-R^{(infty)}(k)}
\end{eqnarray}
\end{lem}

Before the proof of Remark~\ref{lem-limit-G_n}, we give a remark on
$\breve{\vc{G}}$ and $\breve{\vc{R}}(k)$.
\begin{rem}\label{rem-G-matrix}
Consider an M/G/1-type stochastic matrix:
\begin{equation}
\breve{\vc{P}}_{\rm M/G/1}
=
\left(
\begin{array}{ccccc}
\breve{\vc{B}}(0) & 
\breve{\vc{A}}(1) &
\breve{\vc{A}}(2) &
\breve{\vc{A}}(3) &
\cdots
\\
\breve{\vc{A}}(-1)& 
\breve{\vc{A}}(0) &
\breve{\vc{A}}(1) &
\breve{\vc{A}}(2) &
\cdots
\\
\vc{O}&
\breve{\vc{A}}(-1)& 
\breve{\vc{A}}(0) &
\breve{\vc{A}}(1) &
\cdots
\\
\vc{O}					&
\vc{O}					&
\breve{\vc{A}}(-1)& 
\breve{\vc{A}}(0) &
\cdots
\\
\vdots					&
\vdots					&
\ddots					&
\ddots					&
\ddots
\end{array}
\right),
\label{defn-P_{M/G/1}}
\end{equation}
where $\breve{\vc{B}}(0)=\breve{\vc{A}}(-1)+\breve{\vc{A}}(0)$. From
(\ref{eqn-P^{(infty)}}), we have
\[
\breve{\vc{P}}_{\rm M/G/1}
=
\left(
\begin{array}{c|c}
\breve{\vc{B}}(0) &
%
%
\begin{array}{cccc}
\breve{\vc{A}}(1) &
\breve{\vc{A}}(2) &
\breve{\vc{A}}(3) &
\cdots
\end{array}
\\
\hline
%
\begin{array}{c}
\rule{0mm}{5mm}\breve{\vc{A}}(-1)
\\
\vc{O}
\\
\vc{O}
\\
\vdots
\end{array}
%
& ~~\mbox{\LARGE $\breve{\vc{P}}_{\infty}$ }
\end{array}
\right).
\]
Thus (\ref{lim-G_n}) and (\ref{eqn-R^{(infty)}(k)}) imply that
$\breve{\vc{G}}$ and $\breve{\vc{R}}(k)$ are the $G$- and $R$-matrices
of the M/G/1-type stochastic matrix $\breve{\vc{P}}_{\rm M/G/1}$ (see
\cite{Zhao03}). It follows from
(\ref{eqn-breve{A}(-1)})--(\ref{eqn-breve{A}(k)}) and statement (ii)
of Lemma~\ref{overline{breve{A}}(k)>O} that
\[
\breve{\vc{A}}:=\sum_{k=-1}^{\infty}\breve{\vc{A}}(k)
= {\theta \over \mu + \theta}\vc{I}
+ {\mu \over \mu + \theta}\vc{A} > \vc{O},
\]
which shows that $\breve{\vc{A}}$ is an irreducible stochastic matrix
and has the same stationary probability vector $\vc{\pi}$ as
$\vc{A}$. Further combining
(\ref{eqn-breve{A}(-1)})--(\ref{eqn-breve{A}(k)}) with (\ref{eqn-rho})
and $\rho < 1$, we have
\[
\breve{\rho}
:=
\vc{\pi}\sum_{k=0}^{\infty}(k+1)\breve{\vc{A}}(k)\vc{e}
= {\theta \over \mu + \theta} + {\mu \over \mu + \theta}
\vc{\pi}\sum_{k=1}^{\infty}k\vc{A}(k)\vc{e}
= {\rho\mu + \theta \over \mu + \theta}< 1.
\]
Note here that $\breve{\vc{P}}_{\rm M/G/1}$ is irreducible due to
Lemma~\ref{overline{breve{A}}(k)>O}.  These facts imply that the
irreducible stochastic matrix $\breve{\vc{P}}_{\rm M/G/1}$ is positive
recurrent (see, e.g., \cite[Chapter XI, Proposition 3.1]{Asmu03}). In
addition, $\breve{\vc{G}}$ is stochastic (see
\cite[Theorem~2.3.1]{Neut89}) and the spectral radius of
$\breve{\vc{R}}:=\sum_{k=1}^{\infty}\breve{\vc{R}}(k)$ is less than
one (see \cite[Theorem 23]{Zhao03}).
\end{rem}

{\it Proof of Lemma~\ref{lem-limit-G_n}.} Using the dominated
convergence theorem, we take the limit of
(\ref{defn-N_n(0,0)})--(\ref{defn-R_n(k)}) as $n \to \infty$ and
obtain (\ref{lim-G_n})--(\ref{eqn-R^{(infty)}(k)}). Therefore it
remains to prove the positivity of the limiting matrices.

We note that $\breve{\vc{G}}$ is the $G$-matrix of the M/G/1-type
stochastic matrix $\breve{\vc{P}}_{\rm M/G/1}$ (see
Remark~\ref{rem-G-matrix}) and thus it is the unique accumulation
point of the following sequence $\{\breve{\vc{G}}_{\nu}\}$ (see
\cite[Chapter~2]{Neut89}):
\[
\breve{\vc{G}}_0 = \vc{O},\quad
\breve{\vc{G}}_{\nu} =  \breve{\vc{A}}(-1) 
+ \sum_{k=0}^{\infty} \breve{\vc{A}}(k) \left(\breve{\vc{G}}{}_{\nu-1} \right)^{k+1}
~~~\mbox{for $\nu \in \bbN$},
\]
which leads to
\begin{equation}
\breve{\vc{G}}
\ge \breve{\vc{A}}(-1) 
+ \sum_{k=0}^{\infty} \breve{\vc{A}}(k)\left( \breve{\vc{A}}(-1) \right)^{k+1}.
\label{ineq-01}
\end{equation}
 Lemma~\ref{overline{breve{A}}(k)>O} shows that the diagonal elements
 of $\breve{\vc{A}}(-1)$ are all positive and $\sum_{k=0}^{\infty}
 \breve{\vc{A}}(k) > \vc{O}$. Thus from (\ref{ineq-01}) and
 (\ref{eqn-U^{(infty)}(0)}), we have
\[
\breve{\vc{G}} > \vc{O},
\quad \breve{\vc{U}}(0) =
 \sum_{k=0}^{\infty} \breve{\vc{A}}(k)\breve{\vc{G}}{}^k > \vc{O},
\quad ( \vc{I} - \breve{\vc{U}}(0) )^{-1} \ge \breve{\vc{U}}(0) >
 \vc{O}.
\]
Finally, the positivity of $\breve{\vc{R}}(k)$ in
(\ref{eqn-R^{(infty)}(k)}) follows from $\breve{\vc{G}}( \vc{I} -
\breve{\vc{U}}(0) )^{-1} > \vc{O}$ and statement (ii) of
Lemma~\ref{overline{breve{A}}(k)>O}. \qed

\medskip

Lemma~\ref{lem-limit-G_n} and Remark~\ref{rem-G-matrix} show that
$\breve{\vc{G}}$ is an irreducible stochastic matrix.  Thus
$\breve{\vc{G}}$ has the unique and positive stationary probability
vector, which is denoted by $\breve{\vc{g}} > \vc{0}$ hereafter.

Lemma~\ref{lem-g} below shows a relationship between $\breve{\vc{g}}$
and $\vc{\pi}$. We can readily prove this lemma by using
Remark~\ref{rem-G-matrix} and following the proof of Lemma~14 in
\cite{Taki03}. Thus we omit the proof.
\begin{lem}\label{lem-g}
\begin{equation}
\vc{\pi} 
= (1 - \breve{\rho})\breve{\vc{g}} 
\left( \vc{I} - \breve{\vc{U}}(0) \right)^{-1}
\left( \vc{I} - \breve{\vc{R}} \right)^{-1},
\label{eqn-pi-g}
\end{equation}
where $\breve{\vc{R}} = \sum_{k=1}^{\infty}\breve{\vc{R}}(k)$ and
$\breve{\rho} = (\rho\mu + \theta) / (\mu + \theta)$.
\end{lem}

\medskip

Recall that $\vc{q}= (\vc{q}(0),\vc{q}(1),\vc{q}(2),\dots)$ is the
stationary probability vector of $\breve{\vc{P}}$ in
(\ref{eqn-breve{P}}). According to Theorems~2.1 and 2.6 in
\cite{Zhao98}, $\vc{q}= (\vc{q}(0),\vc{q}(1),\vc{q}(2),\dots)$ can be
characterized as follows:
\begin{equation}
\vc{q}(k)
= \sum_{n=0}^{k-1}\vc{q}(n)\breve{\vc{R}}_n(k-n),
\qquad k \in \bbN.
\label{eqn-q(k)}
\end{equation}
Therefore we discuss the asymptotics for $\{\overline{\vc{q}}(k)\}$
through $\{\breve{\vc{R}}_n(k)\}$, which requires some preparations.

\begin{lem}\label{lem-asymp-D(k)}
If Assumption~\ref{assumpt-subexp} holds, then (i) $\vc{c}^{\rm
  D}:=\limsup_{k\to\infty}\ooverline{\vc{D}}(k)\vc{e}/\PP(Y > k)$ is
finite; and (ii) $\lim_{k\to\infty}\overline{\vc{A}}(k)/\PP(Y > k)=
\vc{O}$ and $\lim_{k\to\infty}\overline{\vc{D}}(k)/\PP(Y > k)=
\vc{O}$.
\end{lem}

\proof
From (\ref{defn-A(k)}), we have
\begin{eqnarray*}
\sum_{k=0}^{\infty}z^k \vc{A}(k)
&=& \sum_{m=0}^{\infty}
\int_0^{\infty} \rme^{-\theta x}{(\theta x)^m \over m!} \rmd H(x)
\cdot \left[
\vc{I} + \theta^{-1}\left( \vc{C} + \widehat{\vc{D}}(z) \right) 
\right]^m
\nonumber
\\
&\ge& 
\int_0^{\infty} \rme^{-\theta x}(\theta x) \rmd H(x) \cdot 
\left[
\vc{I} + \theta^{-1}\left( \vc{C} + \widehat{\vc{D}}(z) \right) 
\right],
\end{eqnarray*}
which leads to
\begin{equation}
\vc{A}(k)
\ge \zeta \cdot \vc{D}(k),
\qquad k \in \bbN,
\label{ineqn-A(k)-D(k)}
\end{equation}
where $\zeta = \int_0^{\infty} x\rme^{-\theta x} \rmd H(x) \in
(0,\infty)$ due to $h = \int_0^{\infty} x \rmd H(x)\in
(0,\infty)$. Therefore (\ref{ineqn-A(k)-D(k)}) and
Assumption~\ref{assumpt-subexp} show that statement (i) is
true. Further, Assumption~\ref{assumpt-subexp} implies that
\[
\lim_{k\to\infty}
{\overline{\vc{A}}(k) \over \PP(Y > k) }
\le
\lim_{k\to\infty}
{\overline{\vc{A}}(k)\vc{e}\vc{e}^{\rmt} \over \PP(Y > k) }
= \lim_{k\to\infty}
{
\ooverline{\vc{A}}(k-1)\vc{e}\vc{e}^{\rmt} 
- \ooverline{\vc{A}}(k)\vc{e}\vc{e}^{\rmt} \over \PP(Y > k) 
} 
= \vc{O}.
\]
Combining this and (\ref{ineqn-A(k)-D(k)}) yields
$\lim_{k\to\infty}\overline{\vc{D}}(k)/\PP(Y > k)= \vc{O}$.  \qed

\medskip

It follows from (\ref{eqn-breve{A}-00b}) and (\ref{eqn-breve{A}-03})
that
\begin{equation}
\breve{\vc{A}}_n(k)
= {\min(n,1)\mu \over \mu + \theta} \vc{A}(k+1)
+ {\vc{D} \ast \vc{A}(k+1) \over \max(n,1)(\mu + \theta)},
\quad n \in \bbZ_+,~k\in\bbN,
\label{eqn-breve{A}_n(k)}
\end{equation}
where $\vc{D}(0) = \vc{O}$ is defined for convenience. Using
(\ref{eqn-breve{A}_n(k)}), we show the asymptotics of
$\{\overline{\breve{\vc{A}}}_n(k)\}$ and
$\{\ooverline{\breve{\vc{A}}}_n(k)\}$.
\begin{lem}\label{lem-asymp-A_n(k)}
Suppose that Assumption~\ref{assumpt-subexp} is satisfied, then the
following hold:
\begin{enumerate}
\item For $n \in \bbZ_+$,
\begin{eqnarray}
\lim_{k\to\infty}{\overline{\breve{\vc{A}}}_n(k) \over \PP(Y > k)}
&=& \vc{O},
\label{limit-A_n(k)-01}
\\
\limsup_{k\to\infty}{\ooverline{\breve{\vc{A}}}_n(k)\vc{e} \over \PP(Y > k)}
&\le& {\mu \over \mu + \theta}
\left( \min(n,1) \vc{c}^{\rm A}
+ {\vc{c}^{\rm D} + \vc{D}\vc{c}^{\rm A} \over \max(n,1)\mu}
\right)
\nonumber
\\
&=:& {\mu \over \mu + \theta}\vc{c}_n^{\rm A}, 
\label{limit-A_n(k)-02}
\end{eqnarray}
where $\sup_{n\in\bbZ_+}\vc{c}_n^{\rm A}$ is finite and $\vc{c}_n^{\rm
  A}$ is nonzero for all $n \in \bbN$ (but $\vc{c}_0^{\rm A}$ can be a
zero vector).
\item If $\lim_{k\to\infty}\ooverline{\vc{D}}(k)\vc{e}/\PP(Y>k) =
  \vc{c}^{\rm D}$, then
\begin{equation}
\lim_{k\to\infty}{\ooverline{\breve{\vc{A}}}_n(k)\vc{e} \over \PP(Y > k)}
= {\mu \over \mu + \theta}\vc{c}_n^{\rm A},\qquad n \in \bbZ_+.
\label{lim-A_n(k)}
\end{equation}
\end{enumerate}
\end{lem}

\proof 
From (\ref{eqn-breve{A}_n(k)}), we have
\begin{align}
\overline{\breve{\vc{A}}}_n(k)
&= {\min(n,1)\mu \over \mu + \theta} \,\overline{\vc{A}}(k+1)
+ { \overline{\vc{D} \ast \vc{A}}(k+1) \over \max(n,1)(\mu + \theta)},
& (n,k) &\in \bbZ_+^2,
\label{ineqn-overline{breve{A}}_n(k)}
\\
\ooverline{\breve{\vc{A}}}_n(k)\vc{e}
&= {\min(n,1)\mu \over \mu + \theta} \,\ooverline{\vc{A}}(k+1)\vc{e}
\nonumber
\\
& \qquad {} 
+ {\ooverline{\vc{D}}(k+1)\vc{e} + \overline{\vc{D} \ast \overline{\vc{A}}}(k+1) \vc{e}
\over \max(n,1)(\mu + \theta)},
& (n,k) &\in \bbZ_+^2,
\label{ineqn-ooverline{breve{A}}_n(k)}
\end{align}
where we use $\overline{\vc{D} \ast \vc{A}}(k) =
\overline{\vc{D}}(k)\vc{A} + \vc{D} \ast \overline{\vc{A}}(k)$ in
(\ref{ineqn-ooverline{breve{A}}_n(k)}).  It follows from $Y \in \calS
\subset \calL$, statement (ii) of Lemma~\ref{lem-asymp-D(k)} and
Proposition~\ref{prop-Masu11} that
\begin{equation}
\lim_{k\to\infty}
{\overline{\vc{D} \ast \vc{A}}(k+1) \over \PP(Y > k)} = \vc{O}.
\label{add-eqn-09a}
\end{equation}
Applying (\ref{add-eqn-09a}) and
$\lim_{k\to\infty}\overline{\vc{A}}(k+1)/\PP(Y > k) = \vc{O}$ to
(\ref{ineqn-overline{breve{A}}_n(k)}) yields (\ref{limit-A_n(k)-01}).

Using Assumption~\ref{assumpt-subexp}, Lemma~\ref{lem-asymp-D(k)} and
Proposition~\ref{prop-Masu11}, we obtain
\begin{eqnarray}
\limsup_{k\to\infty}
{\ooverline{\vc{D}}(k+1) \vc{e} + \overline{\vc{D} \ast \overline{\vc{A}}}(k+1)\vc{e}
\over \PP(Y > k)}
&\le& \vc{c}^{\rm D} + \vc{D}\vc{c}^{\rm A}.
\label{add-eqn-09b}
\end{eqnarray}
Further, if $\lim_{k\to\infty}\ooverline{\vc{D}}(k)\vc{e}/\PP(Y>k) =
\vc{c}^{\rm D}$, 
\begin{equation}
\lim_{k\to\infty}
{\ooverline{\vc{D}}(k+1) \vc{e} + \overline{\vc{D} \ast \overline{\vc{A}}}(k+1)\vc{e}
\over \PP(Y > k)}
= \vc{c}^{\rm D} + \vc{D}\vc{c}^{\rm A}.
\label{add-eqn-09c}
\end{equation}
Applying (\ref{add-eqn-09b}) and Assumption~\ref{assumpt-subexp} to
(\ref{ineqn-ooverline{breve{A}}_n(k)}), we have
(\ref{limit-A_n(k)-02}).  Similarly, if
$\lim_{k\to\infty}\ooverline{\vc{D}}(k)\vc{e}/\PP(Y>k) = \vc{c}^{\rm
  D}$, we have (\ref{lim-A_n(k)}), though we use (\ref{add-eqn-09c})
instead of (\ref{add-eqn-09b}). The statement on $\{\vc{c}_n^{\rm
  A}\}$ follows from the definition of $\{\vc{c}_n^{\rm A}\}$ and
$\vc{c}^{\rm A} \ge \vc{0},\neq\vc{0}$.  \qed

\begin{lem}\label{lem-limit-R_n(k)}
If Assumption~\ref{assumpt-subexp} is satisfied, then the
following hold:
\begin{enumerate}
\item The limit
\begin{equation}
\lim_{k\to\infty}{\overline{\breve{\vc{R}}}(k) \over \PP(Y > k)}
= {\mu \over \mu+\theta}\vc{c}^{\rm A} 
\breve{\vc{g}}\left( \vc{I} - \breve{\vc{U}}(0) \right)^{-1} 
=: \vc{C}^{\rm R}
\label{limit-R(k)}
\end{equation}
exists, and $\vc{C}^{\rm R}$ has no zero columns.
\item For $n \in \bbZ_+$,
\begin{eqnarray}
\limsup_{k\to\infty}{\overline{\breve{\vc{R}}}_n(k) \over \PP(Y > k)}
&\le& {\mu \over \mu + \theta}\vc{c}_n^{\rm A}
\breve{\vc{g}}(\vc{I}-\breve{\vc{U}}(0))^{-1}=: \vc{C}_n^{\rm R},
\label{limsup-R_n(k)}
\end{eqnarray}
where $\sup_{n\in\bbZ_+}\vc{C}_n^{\rm R}$ is finite and $\vc{C}_n^{\rm
  R}$ has no zero columns for all $n \in \bbN$ (but $\vc{C}_0^{\rm R}$
can be a zero matrix).
\item If $\lim_{k\to\infty}\ooverline{\vc{D}}(k)\vc{e}/\PP(Y>k) =
  \vc{c}^{\rm D}$, then
\begin{eqnarray}
\lim_{k\to\infty}{\overline{\breve{\vc{R}}}_n(k) \over \PP(Y > k)}
= \vc{C}_n^{\rm R},\qquad n \in \bbZ_+.
\label{lim-R_n(k)}
\end{eqnarray}
\end{enumerate}
\end{lem}

\proof See Appendix~\ref{proof-lem-limit-R_n(k)}. \qed

\begin{lem}\label{lem-ineqn-R_n(k)}
Let $\vc{\varGamma}(k)$ ($k\in\bbZ_+$) denote
\begin{equation}
\vc{\varGamma}(k)
= \sum_{m=0}^{\infty}\vc{D} \ast \vc{A}(k+m+1)
\breve{\vc{G}}{}^m (\vc{I} - \breve{\vc{U}}(0))^{-1},\qquad k \in \bbZ_+.
\label{defn-Gamma(k)}
\end{equation}
The following hold:
\begin{enumerate}
\item For any $\varepsilon > 0$, there exists some
  $n_0:=n_0(\varepsilon) \in \bbN$ such that for all $n \ge n_0$,
\begin{equation}
(1 - \varepsilon)\breve{\vc{R}}(k)
\le \breve{\vc{R}}_n(k)
\le (1 + \varepsilon)
\left\{ \breve{\vc{R}}(k) + \varepsilon  \vc{\varGamma}(k) \right\},
\qquad k \in \bbN.
\label{add-eqn-03}
\end{equation}
\item If Assumption~\ref{assumpt-subexp} holds, then
\begin{equation}
\limsup_{k\to\infty}{\overline{\vc{\varGamma}}(k) \over \PP(Y > k)}
\le (\vc{D}\vc{c}^{\rm A} + \vc{c}^{\rm D})
\breve{\vc{g}}\left( \vc{I} - \breve{\vc{U}}(0) \right)^{-1}
=: \vc{C}^{\rm \Gamma}.
\end{equation}
\end{enumerate}

\end{lem}

\proof See Appendix~\ref{proof-lem-ineqn-R_n(k)}. \qed

We are now ready to prove Lemma~\ref{lem-asymp-q(k)}.

\medskip

\noindent
{\it Proof of Lemma~\ref{lem-asymp-q(k)}}~ For $\varepsilon > 0$, we
fix $n_0:=n_0(\varepsilon)$ for which statement (i) of
Lemma~\ref{lem-ineqn-R_n(k)} holds. We then define
$\vc{s}_{\varepsilon}^+=
(\vc{s}_{\varepsilon}^+(0),\vc{s}_{\varepsilon}^+(1),\vc{s}_{\varepsilon}^+(2),\dots)$
and $\vc{s}_{\varepsilon}^-=
(\vc{s}_{\varepsilon}^-(0),\vc{s}_{\varepsilon}^-(1),\vc{s}_{\varepsilon}^-(2),\dots)$
as follows:
\begin{equation}
\vc{s}_{\varepsilon}^+(0) 
= \vc{s}_{\varepsilon}^-(0) 
= (\vc{q}(0),\vc{q}(1),\dots,\vc{q}(n_0)),
\label{defn-s(0)}
\end{equation}
and for $ k \in \bbN$,
\begin{eqnarray}
\vc{s}_{\varepsilon}^+(k) 
&=& \vc{s}_{\varepsilon}^+(0) \breve{\vc{R}}_{(0,n_0)}(k)
+ (1+\varepsilon)\sum_{n=1}^{k-1}\vc{s}_{\varepsilon}^+(n)
\left( \breve{\vc{R}}(k-n) + \varepsilon\vc{\varGamma}(k-n) \right),
\qquad~~
\label{defn-s^+(k)}
\\
\vc{s}_{\varepsilon}^-(k) 
&=& \vc{s}_{\varepsilon}^-(0) \breve{\vc{R}}_{(0,n_0)}(k)
+ (1-\varepsilon)\sum_{n=1}^{k-1}\vc{s}_{\varepsilon}^-(n)\breve{\vc{R}}(k-n),
\label{defn-s^-(k)}
\end{eqnarray}
where
\begin{equation}
\breve{\vc{R}}_{(0,n_0)}(k) 
=
\left(
\begin{array}{c}
\breve{\vc{R}}_0(k+n_0)
\\
\breve{\vc{R}}_1(k+n_0-1)
\\
\vdots
\\
\breve{\vc{R}}_{n_0}(k)
\end{array}
\right),
\qquad k \in \bbN.
\label{defn-breve{R}_{(0,n_0)}(k)}
\end{equation}

For convenience, let $\breve{\vc{R}}_{(0,n_0)}(0) = \vc{O}$ and
$\breve{\vc{R}}(0)=\vc{O}$. Let
$\breve{\vc{R}}=\sum_{k=0}^{\infty}\breve{\vc{R}}(k)$ and
$\vc{\varGamma}=\sum_{k=0}^{\infty}\vc{\varGamma}(k)$. Recall here
that the spectral radius of $\breve{\vc{R}}$ is less than one (see
Remark~\ref{rem-G-matrix}) and thus so is that of
$(1-\varepsilon)\breve{\vc{R}}$. Further for any sufficiently small
$\varepsilon > 0$, the spectral radius of
$(1+\varepsilon)(\breve{\vc{R}} + \varepsilon \vc{\varGamma})$ is less
than one (see, e.g., Theorem~8.1.18 in \cite{Horn90}). We fix
$\varepsilon > 0$ to be such a small value.

Following the proof of Theorem~1 in \cite{Taki04}, we can readily show
that
\begin{equation}
\vc{s}_{\varepsilon}^+(k) 
= \vc{s}_{\varepsilon}^+(0) \breve{\vc{R}}_{(0,n_0)} 
\ast \sum_{m=0}^{\infty} (1+\varepsilon)^m 
(\breve{\vc{R}} + \varepsilon \vc{\varGamma})^{\ast m}(k),
\qquad k \in \bbN,
\label{eqn-s(k)}
\end{equation}
where $\{(\breve{\vc{R}} + \varepsilon \vc{\varGamma})^{\ast
  m}(k);k\in\bbZ_+\}$ is the $m$-fold convolution of
$\{\breve{\vc{R}}(k) + \varepsilon \vc{\varGamma}(k);k\in\bbZ_+\}$
itself. It follows from statement (i) of Lemma~\ref{lem-limit-R_n(k)},
statement (ii) of Lemma~\ref{lem-ineqn-R_n(k)} and
Proposition~\ref{prop-Jele98} that
\begin{eqnarray}
\lefteqn{
\limsup_{k\to\infty} 
\sum_{m=0}^{\infty} {(1+\varepsilon)^m 
\overline{(\breve{\vc{R}} + \varepsilon \vc{\varGamma})^{\ast m}}(k) \over \PP(Y>k)}
}
\quad &&
\nonumber
\\
&\le& \left\{\vc{I} - (1+\varepsilon)
(\breve{\vc{R}} + \varepsilon\vc{\varGamma})\right\}^{-1} 
(1+\varepsilon)(\vc{C}^{\rm R} + \varepsilon \vc{C}^{\rm \Gamma})
\nonumber
\\
&& {} \times 
\left\{\vc{I} - (1+\varepsilon)(\breve{\vc{R}} + \varepsilon\vc{\varGamma})\right\}^{-1}.
\label{eqn-05}
\end{eqnarray}
Further statement (ii) of Lemma~\ref{lem-limit-R_n(k)} yields
\begin{equation}
\limsup_{k\to\infty}
{\overline{\breve{\vc{R}}}_{(0,n_0)}(k) \over \PP(Y > k)}
\le 
\left(
\begin{array}{c}
\vc{C}_0^{\rm R}
\\
\vc{C}_1^{\rm R}
\\
\vdots
\\
\vc{C}_{n_0}^{\rm R} 
\end{array}
\right)
=: \vc{C}_{(0,n_0)}^{\rm R} \neq \vc{O}.
\label{eqn-04}
\end{equation}
Applying Proposition~\ref{prop-Masu11} to (\ref{eqn-s(k)}) and using
(\ref{eqn-05}) and (\ref{eqn-04}), we obtain
\begin{eqnarray}
\lefteqn{
\limsup_{k\to\infty}
{\overline{\vc{s}}_{\varepsilon}^+(k) \over \PP(Y > k)}
}
\quad &&
\nonumber
\\
&\le& \vc{s}_{\varepsilon}^+(0)
\vc{C}_{(0,n_0)}^{\rm R} 
\left\{\vc{I} - (1+\varepsilon)(\breve{\vc{R}} + \varepsilon\vc{\varGamma})\right\}^{-1}
\nonumber
\\
&& {} + \vc{s}_{\varepsilon}^+(0)
\sum_{n=1}^{\infty}\breve{\vc{R}}_{(0,n_0)}(n) 
\left\{
\vc{I} - (1+\varepsilon)(\breve{\vc{R}} + \varepsilon\vc{\varGamma})
\right\}^{-1}
\nonumber
\\
&& \qquad\qquad {} \times 
(1+\varepsilon)(\vc{C}^{\rm R} + \varepsilon \vc{C}^{\rm \Gamma})
\left\{
\vc{I} - (1+\varepsilon)(\breve{\vc{R}} + \varepsilon\vc{\varGamma})
\right\}^{-1}.
\label{eqn-06}
\end{eqnarray}
Recall here that $n_0 \to \infty$ as $\varepsilon \downarrow 0$ (see
Lemmas~\ref{lem-limit-G_n} and \ref{lem-ineqn-R_n(k)}). Recall also
that $\sup_{n\in\bbZ_+} \sum_{k=1}^{\infty}\breve{\vc{R}}_n(k)\vc{e}$
is finite (see Lemma~\ref{lem-RG-bounds}). Thus using
(\ref{defn-s(0)}), (\ref{defn-breve{R}_{(0,n_0)}(k)}) and the
dominated convergence theorem, we have
\begin{eqnarray*}
\lim_{\varepsilon\downarrow0}
\vc{s}_{\varepsilon}^+(0)
\sum_{n=1}^{\infty}\breve{\vc{R}}_{(0,n_0)}(n)
&=& \lim_{n_0\to\infty} 
\sum_{n=1}^{\infty} 
\sum_{l=0}^{n_0} \vc{q}(l)\breve{\vc{R}}_l(n+n_0-l)
\nonumber
\\
&=& \lim_{n_0\to\infty} 
\sum_{l=0}^{n_0} \vc{q}(l)\overline{\breve{\vc{R}}}_l(n_0-l)
\nonumber
\\
&=& 
\sum_{l=0}^{\infty} \vc{q}(l)
\lim_{n_0\to\infty} \overline{\breve{\vc{R}}}_l(n_0-l)
= \vc{0}.
\end{eqnarray*}
Therefore letting $\varepsilon \downarrow 0$ in (\ref{eqn-06}) and
using (\ref{defn-s(0)}) and (\ref{defn-breve{R}_{(0,n_0)}(k)}) yield
\begin{eqnarray}
\lim_{\varepsilon\downarrow0}
\limsup_{k\to\infty}
{\overline{\vc{s}}_{\varepsilon}^+(k) \over \PP(Y > k)}
&\le&
\lim_{\varepsilon\downarrow0}
\vc{s}_{\varepsilon}^+(0)\vc{C}_{(0,n_0)}^{\rm R}
\left( \vc{I} - \breve{\vc{R}} \right)^{-1}
\nonumber
\\
&=& 
\sum_{n=0}^{\infty}\vc{q}(n)\vc{C}_n^{\rm R}
\left( \vc{I} - \breve{\vc{R}} \right)^{-1}.
\label{eqn-07}
\end{eqnarray}
It also follows from (\ref{eqn-pi-g}) and the definition
of $\vc{C}_n^{\rm R}$ (see (\ref{limsup-R_n(k)})) that
\begin{eqnarray*}
\vc{C}_n^{\rm R}\left( \vc{I} -\breve{\vc{R}} \right)^{-1}
= {\mu \over \mu+\theta} {\vc{c}_n^{\rm A}\vc{\pi} \over 1 - \breve{\rho}}.
\end{eqnarray*}
Substituting this equation into (\ref{eqn-07}), we obtain
\begin{eqnarray}
\lim_{\varepsilon\downarrow0}
\limsup_{k\to\infty}
{\overline{\vc{s}}_{\varepsilon}^+(k) \over \PP(Y > k)}
&\le& \sum_{n=0}^{\infty}\vc{q}(n) \vc{c}_n^{\rm A}
\cdot {\mu \over \mu+\theta} {\vc{\pi} \over 1 - \breve{\rho}}.
\label{eqn-08a}
\end{eqnarray}

It is proved later that
\begin{equation}
\vc{s}_{\varepsilon}^-(k) 
\le \vc{q}(k+n_0) \le \vc{s}_{\varepsilon}^+(k),
\qquad k\in\bbN.
\label{ineqn-s(k)}
\end{equation}
Combining (\ref{ineqn-s(k)}) with
(\ref{eqn-08a}) and using $Y \in \calS \subset \calL$, we obtain
\begin{equation}
\limsup_{k\to\infty}
{ \overline{\vc{q}}(k) \over \PP(Y > k) }
\le \sum_{n=0}^{\infty}\vc{q}(n)\vc{c}_n^{\rm A}
\cdot {\mu \over \mu+\theta}{\vc{\pi} \over 1 - \breve{\rho}}.
\label{limsup-overline{q}(k)}
\end{equation}
Note here that $\vc{q}(n) > \vc{0}$ for all $n \in \bbZ_+$ (see
(\ref{defn-q(k)})) and that $\sup_{n\in\bbZ_+}\vc{c}_n^{\rm A}$ is
finite and $\vc{c}_n^{\rm A} \ge \vc{0},\neq\vc{0}$ for all $n\in\bbN$
(see statement (i) of Lemma~\ref{lem-asymp-A_n(k)}). As a result,
\[
0< \sum_{n=0}^{\infty}\vc{q}(n)\vc{c}_n^{\rm A} < \infty,
\]
which completes the proof of statement (i).

Next we prove statement (ii) under the condition that
$\lim_{k\to\infty}\ooverline{\vc{D}}(k)\vc{e} = \vc{c}^{\rm D}$. As
with (\ref{eqn-s(k)}), the following equation holds:
\begin{equation}
\vc{s}_{\varepsilon}^-(k) 
= \vc{s}_{\varepsilon}^-(0) \breve{\vc{R}}_{(0,n_0)} 
\ast \sum_{m=0}^{\infty} (1-\varepsilon)^m \breve{\vc{R}}{}^{\ast m}(k),
\qquad k \in \bbN.
\label{eqn-s^-(k)}
\end{equation}
It follows from statements (i) and (iii) of
Lemma~\ref{lem-limit-R_n(k)} and Proposition~\ref{prop-Jele98} that
\begin{eqnarray}
\lim_{k\to\infty}
{\overline{\breve{\vc{R}}}_{(0,n_0)}(k) \over \PP(Y > k)}
&=&  \vc{C}_{(0,n_0)}^{\rm R},
\label{eqn-04b}
\\
\lim_{k\to\infty}
\sum_{m=0}^{\infty}
{
 (1-\varepsilon)^m \overline{\breve{\vc{R}}{}^{\ast m}}(k)
\over \PP(Y > k)}
&=& \left\{\vc{I} - (1-\varepsilon)
\breve{\vc{R}}\right\}^{-1} 
(1-\varepsilon) \vc{C}^{\rm R}
\nonumber
\\
&& {} \times 
\left\{\vc{I} - (1-\varepsilon)\breve{\vc{R}}\right\}^{-1}.
\label{eqn-04c}
\end{eqnarray}
Using (\ref{eqn-s^-(k)})--(\ref{eqn-04c}) and following the proof of
statement (i), we can show that
\begin{eqnarray*}
\liminf_{k\to\infty}
{ \overline{\vc{q}}(k) \over \PP(Y > k) }
\ge \lim_{\varepsilon\downarrow0}
\lim_{k\to\infty}
{\overline{\vc{s}}_{\varepsilon}^-(k) \over \PP(Y > k)}
= \sum_{n=0}^{\infty}\vc{q}(n) \vc{c}_n^{\rm A}
\cdot {\mu \over \mu+\theta} {\vc{\pi} \over 1 - \breve{\rho}}.
\end{eqnarray*}
This inequality and (\ref{limsup-overline{q}(k)}) show that statement
(ii) holds.

Finally, we prove (\ref{ineqn-s(k)}) by induction.  From
(\ref{eqn-q(k)}), (\ref{defn-s^+(k)}) and (\ref{defn-s^-(k)}), we have
\begin{equation}
\vc{s}_{\varepsilon}^+(1) 
= \vc{s}_{\varepsilon}^-(1) 
= \sum_{n=0}^{n_0}\vc{q}(n)\breve{\vc{R}}_n(n_0+1-n)
= \vc{q}(n_0+1),
\label{eqn-s_{delta}(1)}
\end{equation}
which shows that (\ref{ineqn-s(k)}) holds for $k=1$. Suppose that
(\ref{ineqn-s(k)}) holds for some $k=k_{\ast} \in \bbN$. Substituting
this inductive assumption and the right inequality in
(\ref{add-eqn-03}) into (\ref{defn-s^+(k)}) with $k=k_{\ast}+1$ yields
\begin{eqnarray*}
\vc{s}_{\varepsilon}^+(k_{\ast}+1) 
&\ge& \vc{s}_{\varepsilon}^+(0) \breve{\vc{R}}_{(0,n_0)}(k_{\ast}+1)
\nonumber
\\
&& {} 
+ (1+\varepsilon)\sum_{n=1}^{k_{\ast}}\vc{q}(n+n_0)
\left(\breve{\vc{R}}(k_{\ast}+1-n) + \varepsilon\vc{\varGamma}(k_{\ast}+1-n)
\right)
\nonumber
\\
&\ge& \vc{s}_{\varepsilon}^+(0) \breve{\vc{R}}_{(0,n_0)}(k_{\ast}+1) 
+ \sum_{n=1}^{k_{\ast}}\vc{q}(n+n_0)
\breve{\vc{R}}_{n+n_0}(k_{\ast}+1-n)
\nonumber
\\
&=& \sum_{n=0}^{n_0}\vc{q}(n)\breve{\vc{R}}_n(k_{\ast}+1+n_0-n)
+ \sum_{n=1}^{k_{\ast}}\vc{q}(n+n_0)
\breve{\vc{R}}_{n+n_0}(k_{\ast}+1-n)
\nonumber
\\
&=& \sum_{n=0}^{k_{\ast}+n_0}\vc{q}(n)\breve{\vc{R}}_n(k_{\ast}+1+n_0-n)
\nonumber
\\
&=& \vc{q}(k_{\ast}+1+n_0),
\end{eqnarray*}
where the last equality is due to (\ref{eqn-q(k)}). As a result, the
right inequality in (\ref{ineqn-s(k)}) has been proved. The left one
is proved in a similar way.

\appendix

\section{Proofs}

\subsection{Proof of Lemma~\ref{lem-limit-R_n(k)}}\label{proof-lem-limit-R_n(k)} 

It follows from Assumption~\ref{assumpt-subexp} and
(\ref{eqn-breve{A}(k)}) that
\[
\lim_{k\to\infty}{\ooverline{\breve{\vc{A}}}(k) \over \PP(Y > k)}\vc{e} 
= {\mu \over \mu + \theta} \vc{c}^{\rm A}.
\]
Using this equation and proceeding as in the proof of Lemma~3.2 in
\cite{Masu13-ANNOR}, we can show that 
\[
\lim_{k\to\infty}{\overline{\breve{\vc{R}}}(k) \over \PP(Y > k)}
= {\mu \over \mu+\theta}\vc{c}^{\rm A} 
{\vc{\pi}(\vc{I} - \breve{\vc{R}}) \over 1 - \breve{\rho}},
\]
from which and Lemma~\ref{lem-g} it follows that the limit in
(\ref{limit-R(k)}) exists. It also follows from $\breve{\vc{g}} >
\vc{0}$, $( \vc{I} - \breve{\vc{U}}(0))^{-1} > \vc{O}$ (see
Lemma~\ref{lem-limit-G_n}) and $\vc{c}^{\rm A} \ge \vc{0},\neq\vc{0}$
(see Assumption~\ref{assumpt-subexp}) that $\vc{C}^{\rm R}$ has no
zero columns. Thus statement (i) holds.

Similarly we can prove statements (ii) and (iii), though we need
additional steps. For completeness, we provide the proof of statements
(ii) and (iii).

Lemma~\ref{lem-limit-G_n} implies that
\begin{equation}
\lim_{k\to\infty}\left(
\prod_{l=n+k+m}^{n+k+1}\breve{\vc{G}}_l
\right) \breve{\vc{N}}_{n+k}(0,0) 
=\breve{\vc{G}}{}^m \left( \vc{I} - \breve{\vc{U}}(0) \right)^{-1}
\quad\mbox{uniformly over $m,n \in \bbZ_+$}.
\label{add-enq-07}
\end{equation}
Fix $\varepsilon > 0$ arbitrarily, which is independent of $n$. It
then follows from (\ref{add-enq-07}) and (\ref{defn-R_n(k)}) that for all
sufficiently large $k$,
\begin{align}
\overline{\breve{\vc{R}}}_n(k) 
&\le (1+\varepsilon) \sum_{m=0}^{\infty} \overline{\breve{\vc{A}}}_n(k+m)
\breve{\vc{G}}{}^m \left( \vc{I} - \breve{\vc{U}}(0) \right)^{-1},
& n &\in \bbZ_+,
\label{ineqn-05a}
\\
\overline{\breve{\vc{R}}}_n(k) 
&\ge (1-\varepsilon) \sum_{m=0}^{\infty} \overline{\breve{\vc{A}}}_n(k+m)
\breve{\vc{G}}{}^m \left( \vc{I} - \breve{\vc{U}}(0) \right)^{-1},
& n &\in \bbZ_+.
\label{ineqn-05b}
\end{align}
Recall here that $\breve{\vc{G}}$ is a positive stochastic matrix with
stationary probability vector $\breve{\vc{g}}$ and thus
$\lim_{m\to\infty}\breve{\vc{G}}{}^m = \vc{e}\breve{\vc{g}}$ (see,
e.g., \cite[Theorem 8.2.8]{Horn90}). Therefore there exists some $m_1
:=m_1(\varepsilon) \in \bbN$ such that for all $m>m_1$,
\begin{equation}
(1-\varepsilon)\vc{e}\breve{\vc{g}}
\le
\breve{\vc{G}}{}^m
\le (1+\varepsilon)\vc{e}\breve{\vc{g}}.
\label{ineqn-breve{G}^m}
\end{equation}
Substituting the right inequality in (\ref{ineqn-breve{G}^m}) into
(\ref{ineqn-05a}), we have for all sufficiently large $k$,
\begin{eqnarray}
\overline{\breve{\vc{R}}}_n(k) 
&\le& (1+\varepsilon) \sum_{m=0}^{m_1} \overline{\breve{\vc{A}}}_n(k+m)
\breve{\vc{G}}{}^m \left( \vc{I} - \breve{\vc{U}}(0) \right)^{-1}
\nonumber
\\
&& {} 
+ (1+\varepsilon)^2
\sum_{m=m_1+1}^{\infty} \overline{\breve{\vc{A}}}_n(k+m)
\vc{e}\breve{\vc{g}} \left( \vc{I} - \breve{\vc{U}}(0) \right)^{-1}
\nonumber
\\
&\le& (1+\varepsilon) \sum_{m=0}^{m_1} \overline{\breve{\vc{A}}}_n(k+m)
\breve{\vc{G}}{}^m \left( \vc{I} - \breve{\vc{U}}(0) \right)^{-1}
\nonumber
\\
&& {} 
+ (1+\varepsilon)^2 \ooverline{\breve{\vc{A}}}_n(k+m_1)
\vc{e}\breve{\vc{g}} \left( \vc{I} - \breve{\vc{U}}(0) \right)^{-1}.
\label{ineqn-breve{R}_n(k)}
\end{eqnarray}
It also follows from (\ref{limit-A_n(k)-01}) and $Y \in \calS
\subset \calL$ that for any fixed $m \in \bbN$ and $n \in \bbZ_+$,
\begin{equation}
\lim_{k\to\infty}{\overline{\breve{\vc{A}}}_n(k+m) \over \PP(Y > k)}
= \lim_{k\to\infty}{\overline{\breve{\vc{A}}}_n(k+m) \over \PP(Y > k+m)}
{\PP(Y > k+m) \over \PP(Y > k)}
= \vc{O}.
\label{lim-overline{A}(k)}
\end{equation}
Applying (\ref{lim-overline{A}(k)}), (\ref{limit-A_n(k)-02}) and
Proposition~\ref{prop-Masu11} to (\ref{ineqn-breve{R}_n(k)}), we
obtain
\begin{eqnarray}
\limsup_{k\to\infty}
{\overline{\breve{\vc{R}}}_n(k) \over \PP(Y>k) }
&\le&  (1+\varepsilon)^2 {\mu \over \mu + \theta}\vc{c}_n^{\rm A}
\breve{\vc{g}} \left( \vc{I} - \breve{\vc{U}}(0) \right)^{-1}
= (1+\varepsilon)^2 \vc{C}_n^{\rm R}.
\label{add-eqn-06a}
\end{eqnarray}
Letting $\varepsilon \downarrow 0$ in (\ref{add-eqn-06a}) yields
(\ref{limsup-R_n(k)}). In addition, since
$\sup_{n\in\bbZ_+}\vc{c}_n^{\rm A}$ is finite (see
Lemma~\ref{lem-asymp-A_n(k)}), so is $\sup_{n\in\bbZ_+}\vc{C}_n^{\rm
  R}$. In addition, $\vc{C}_n^{\rm R}$ ($\forall n\in\bbN$) has no
zero columns because $\breve{\vc{g}} > \vc{0}$, $( \vc{I} -
\breve{\vc{U}}(0) )^{-1}> \vc{O}$ and $\vc{c}_n^{\rm A} \ge
\vc{0},\neq\vc{0}$ for $n\in\bbN$.

Finally, we assume that $\lim_{k\to\infty} \ooverline{\vc{D}}(k)\vc{e}
/ \PP(Y>k) = \vc{c}^{\rm D}$. Using (\ref{ineqn-05b}) and the left
inequality in (\ref{ineqn-breve{G}^m}) (and following the proof of
(\ref{limsup-R_n(k)})), we can show that
\[
\liminf_{k\to\infty}
{\overline{\breve{\vc{R}}}_n(k) \over \PP(Y>k) }
\ge \vc{C}_n^{\rm R}.
\]
Combining this and (\ref{limsup-R_n(k)}), we have (\ref{lim-R_n(k)}).

\subsection{Proof of Lemma~\ref{lem-ineqn-R_n(k)}}\label{proof-lem-ineqn-R_n(k)}
We estimate $\breve{\vc{R}}_n(k)$ in (\ref{defn-R_n(k)}) as a function
of $n$.  Similarly to (\ref{add-enq-07}), Lemma~\ref{lem-limit-G_n}
implies that
\[
\lim_{n\to\infty}\left(
\prod_{l=n+k+m}^{n+k+1}\breve{\vc{G}}_l
\right) \breve{\vc{N}}_{n+k}(0,0) 
=\breve{\vc{G}}{}^m \left( \vc{I} - \breve{\vc{U}}(0) \right)^{-1},
\]
where the convergence is uniform over $(k,m) \in \bbN \times
\bbZ_+$. Thus for any $\varepsilon > 0$, there exists some
$n':=n'(\varepsilon) \in \bbN$ such that for all $n \ge n'$,
\begin{align}
\breve{\vc{R}}_n(k) 
&\le (1+\varepsilon) \sum_{m=0}^{\infty} \breve{\vc{A}}_n(k+m)
\breve{\vc{G}}{}^m \left( \vc{I} - \breve{\vc{U}}(0) \right)^{-1},
& k &\in \bbN,
\label{ineqn-08a}
\\
\breve{\vc{R}}_n(k) 
&\ge (1-\varepsilon) \sum_{m=0}^{\infty} \breve{\vc{A}}_n(k+m)
\breve{\vc{G}}{}^m \left( \vc{I} - \breve{\vc{U}}(0) \right)^{-1},
& k &\in \bbN.
\label{ineqn-08b}
\end{align}

It follows from (\ref{eqn-breve{A}_n(k)}) and (\ref{eqn-breve{A}(k)})
that for all $n \ge \lceil 1/\{\varepsilon(\mu + \theta)\} \rceil$,
\begin{equation}
\breve{\vc{A}}_n(k) \le \breve{\vc{A}}(k)
+ \varepsilon \vc{D} \ast \vc{A}(k+1),\qquad k \in \bbN,
\label{right-ineqn-breve{A}_n(k)}
\end{equation}
and that for all $n \in \bbN$,
\begin{equation}
\breve{\vc{A}}_n(k) \ge \breve{\vc{A}}(k),
\qquad k \in \bbN.
\label{left-ineqn-breve{A}_n(k)}
\end{equation}
Substituting (\ref{right-ineqn-breve{A}_n(k)}) and
(\ref{left-ineqn-breve{A}_n(k)}) into (\ref{ineqn-08a}) and
(\ref{ineqn-08b}) respectively and using (\ref{eqn-R^{(infty)}(k)})
and (\ref{defn-Gamma(k)}), we obtain for all $n \ge n_0
:=n_0(\varepsilon) = \max(n',\lceil 1/\{\varepsilon(\mu + \theta)\}
\rceil)$,
\begin{align*}
\breve{\vc{R}}_n(k) 
&\le (1+\varepsilon) 
\left( \breve{\vc{R}}(k) 
+ \varepsilon \vc{\varGamma}(k)
\right),
& k &\in \bbN,
\\
\breve{\vc{R}}_n(k) 
&\ge (1-\varepsilon)\breve{\vc{R}}(k),
& k &\in \bbN,
\end{align*}
which show that statement (i) holds.

As for statement (ii), we can prove this by using
(\ref{ineqn-breve{G}^m}), (\ref{add-eqn-09a}) and (\ref{add-eqn-09b})
and following the proof of Lemma~3.2 in \cite{Masu13-ANNOR}. The proof
of statement (ii) is also similar to that of
Lemma~\ref{lem-limit-R_n(k)} (see
Appendix~\ref{proof-lem-limit-R_n(k)}). Therefore we omit the details.

\section{Convolution of Matrix Sequences with Subexponential Tails}\label{appen-subexp}

The following are basic asymptotic results on the convolution of
matrix sequences associated with subexponential tails.
\begin{prop}\label{prop-Jele98}
Suppose that $\{\vc{M}(k);k\in\bbZ_+\}$ is a sequence of nonnegative
square matrices such that $\sum_{n=0}^{\infty}\vc{M}^n = (\vc{I} -
\vc{M})^{-1} < \infty$.
\begin{enumerate}
\item If there exists some $U \in \calS$ such that
\[
\limsup_{k\to\infty}{\overline{\vc{M}}(k) \over \PP(U>k)} 
\le \widetilde{\vc{M}},
\]
then
\[
\limsup_{k\to\infty} {\overline{\sum_{n=0}^{\infty}\vc{M}^{\ast n}}(k) 
\over \PP(U>k)} 
\le (\vc{I} - \vc{M})^{-1}\widetilde{\vc{M}}
(\vc{I} - \vc{M})^{-1}.
\]
\item Replacing ``$\limsup$" and ``$\le$" by ``$\liminf$" and
  ``$\ge$", respectively, in statement (i), we have a true statement.
\item Replacing ``$\limsup$" and ``$\le$" by ``$\lim$" and ``$=$",
  respectively, in statement (i), we have a true statement.
\end{enumerate}

\end{prop}

\begin{prop}\label{prop-Masu11}
Suppose that $\{\vc{M}(k);k\in\bbZ_+\}$ and $\{\vc{N}(k);k\in\bbZ_+\}$
are finite-dimensional nonnegative matrix sequences such that their
convolution is well-defined. Further suppose that $\vc{M}:=
\sum_{k=0}^{\infty}\vc{M}(k)$ and $\vc{N}:=
\sum_{k=0}^{\infty}\vc{N}(k)$ are finite. Under these conditions, the
following hold:
\begin{enumerate}
\item If there exists some $U
\in \calS$ such that
\[
\limsup_{k\to\infty}{\overline{\vc{M}}(k) \over \PP(U > k)} 
\le \widetilde{\vc{M}},
\qquad
\limsup_{k\to\infty}{\overline{\vc{N}}(k) \over \PP(U > k)} 
\le \widetilde{\vc{N}},
\]
then
\[
\limsup_{k\to\infty} 
{\overline{\vc{M} \ast \vc{N}}(k) \over \PP(U > k)} 
\le \widetilde{\vc{M}} \vc{N} 
+ \vc{M} \widetilde{\vc{N}}.
\]
\item Replacing ``$\limsup$" and ``$\le$" by ``$\liminf$" and
  ``$\ge$", respectively, in statement (i), we have a true statement.
\item Replacing ``$\limsup$" and ``$\le$" by ``$\lim$" and ``$=$",
  respectively, in statement (i), we have a true statement.
\end{enumerate}
 
\end{prop}

\noindent
{\it Proof of Propositions~\ref{prop-Jele98} and \ref{prop-Masu11}.~}
The first statements (on the limit superiors) of
Propositions~\ref{prop-Jele98} and \ref{prop-Masu11} are presented in
Lemma~A.12 in \cite{Masu13}. Following the proof of the lemma, we can
readily prove the second statements (on the limit inferiors) of the
two propositions. The third statements are immediate from the first
and second ones, and they also presented in Lemma 6 in \cite{Jele98}
and Proposition~A.3 in \cite{Masu11}.  \qed

\section*{Acknowledgments}
The author thanks Daigo Yoshikawa for pointing out some errors in an
earlier version of this paper. The author is also grateful to
anonymous referees for their invaluable comments and suggestions on
improving the presentation of this paper. Research of the author was
supported in part by Grant-in-Aid for Young Scientists (B) of Japan
Society for the Promotion of Science under Grant No. 24710165.

%
\bibliographystyle{spmpsci}       
%
%

\end{document}